\documentclass[11pt]{article}

\usepackage{graphicx}
\usepackage{color,verbatim}

\usepackage{algorithm}
\usepackage{algorithmic}

\newtheorem{thm}{Theorem}[section]

\newtheorem{pro}[thm]{Proposition}
\newtheorem{lem}[thm]{Lemma}
\newtheorem{cor}[thm]{Corollary}
\newtheorem{alg}[thm]{Algorithm}
\newtheorem{ass}[thm]{Assumption}
\newtheorem{defi}[thm]{Definition}
\newtheorem{example}[thm]{Example}

\newcommand{\sect}[1]{
        \par
        \stepcounter{section}
        \settowidth{\hangindent}{\large\bf\thesection.~}
        \hangafter=1
        \bigskip\bigskip\noindent
        {\large\bf\hbox{\thesection.~}#1}\par
        \nopagebreak
        \medskip
        \renewcommand{\theequation}{\thesection.\arabic{equation}}
        \setcounter{equation}{0}
        \setcounter{subsection}{0}
}
\renewcommand{\subsection}[1]{
        \stepcounter{subsection}
        \noindent
        {\bf\hbox{\thesubsection.~}#1}
        \nobreak
}
\renewcommand{\subsubsection}[1]{
        \stepcounter{subsubsection}
        \noindent
        {\bf\hbox{\thesubsubsection.~}#1}
        \nobreak
}

\newcommand{\nn}{\nonumber}
\newcommand{\be}{\begin{equation}}
\newcommand{\ee}{\end{equation}}
\newcommand{\ba}{\begin{array}}
\newcommand{\ea}{\end{array}}
\newcommand{\bea}{\begin{eqnarray}}
\newcommand{\eea}{\end{eqnarray}}
\newcommand{\bal}{\begin{alg}}
\newcommand{\eal}{\end{alg}}
\newcommand{\ble}{\begin{lem}}
\newcommand{\ele}{\end{lem}}
\newcommand{\bco}{\begin{cor}}
\newcommand{\eco}{\end{cor}}
\newcommand{\bde}{\begin{defi}}
\newcommand{\ede}{\end{defi}}
\newcommand{\bth}{\begin{thm}}
\newcommand{\eth}{\end{thm}}
\newcommand{\bpr}{\begin{pro}}
\newcommand{\epr}{\end{pro}}
\newcommand{\bas}{\begin{ass}}
\newcommand{\eas}{\end{ass}}
\newcommand{\bex}{\begin{example}}
\newcommand{\eex}{\end{example}}
\newcommand{\bfig}{\begin{figure}}
\newcommand{\efig}{\end{figure}}
\newcommand{\reff}[1]{(\ref{#1})}
\newcommand{\refa}[1]{Assumption\ \ref{#1}}
\newcommand{\refl}[1]{Lemma\ \ref{#1}}
\newcommand{\reft}[1]{Theorem\ \ref{#1}}

\newcommand{\refal}[1]{Algorithm\ \ref{#1}}

\newtheorem{prc}[thm]{Procedure}
\newcommand{\bprc}{\begin{prc}}
\newcommand{\eprc}{\end{prc}}

\def\dd{&\!\!\!\!}
\def\eop{\hfill\vbox{\hrule height0.6pt\hbox{\vrule height1.3ex
width0.6pt\hskip1.2ex\vrule width0.6pt}\hrule height0.6pt}}
\def\prf{\noindent {\sl Proof.} \rm}
\def\alglist{
\begin{list}{Step 1}
{\setlength{\leftmargin}{0.5 in}\setlength{\labelwidth}{0.7 in}}
}
\def\eli{\end{list}}
\def\bdes{\begin{description}}
\def\edes{\end{description}}

\def\na{\nabla}

\def\hf{\frac{1}{2}}

\def\alp{\alpha}

\def\st{\hbox{s.t.}}
\def\diag{\hbox{diag}\,}

\begin{document}
\pagenumbering{arabic}
\begin{titlepage}\setcounter{page}{0}

\title{A novel augmented Lagrangian method of multipliers
for optimization with general inequality constraints
}
\author{
Xin-Wei Liu,\thanks{Institute of Mathematics, Hebei University of Technology, Tianjin 300401, China. E-mail:
mathlxw@hebut.edu.cn. The research is supported by the NSFC grants (nos. 12071108 and 11671116).}
\
Yu-Hong Dai,\thanks{Academy of Mathematics and Systems Science, Chinese Academy of Sciences, Beijing 100190, China \& School of Mathematical Sciences, University of Chinese Academy of Sciences, Beijing 100049, China. This author is supported by the NSFC grants (nos. 12021001, 11991021, 11991020 and 11971372),
 the Strategic Priority Research Program of Chinese Academy of Sciences (no. XDA27000000) and
 Beijing Academy of Artificial Intelligence.}
\
Ya-Kui Huang,\thanks{Institute of Mathematics, Hebei University of Technology, Tianjin 300401, China.
This author is supported by the NSFC grant (no. 11701137) and HNSF grant (no. A2021202010).}
\
and Jie Sun\thanks{Institute of Mathematics, Hebei University of Technology, Tianjin 300401, China, and School of Business, National University of Singapore, Singapore 119245, Singapore.}
}
\maketitle

\noindent\underline{\hspace*{6.3in}}
\par

\vskip 10 true pt \noindent{\small{\bf Abstract.}
We introduce a twice differentiable augmented Lagrangian for nonlinear optimization with general inequality constraints and show that a strict local minimizer of the original problem is an approximate strict local solution of the augmented Lagrangian. A novel augmented Lagrangian method of multipliers (ALM) is then presented. Our method is originated from a generalization of the Hetenes-Powell augmented Lagrangian, and is a combination of the augmented Lagrangian and the interior-point technique. It shares a similar algorithmic framework with existing ALMs for optimization with inequality constraints, but it can use the second derivatives and does not depend on projections on the set of inequality constraints. In each iteration, our method solves a twice continuously differentiable unconstrained optimization subproblem on primal variables. The dual iterates, penalty and smoothing parameters are updated adaptively. The global and local convergence are analyzed. Without assuming any constraint qualification, it is proved that the proposed method has strong global convergence. The method may converge to either a Kurash-Kuhn-Tucker (KKT) point or a singular stationary point when the converging point is a minimizer. It may also converge to an infeasible stationary point of nonlinear program when the problem is infeasible. Furthermore, our method is capable of rapidly detecting the possible infeasibility of the solved problem. Under suitable conditions, it is locally linearly convergent to the KKT point, which is consistent with ALMs for optimization with equality constraints. The preliminary numerical experiments on some small benchmark test problems demonstrate our theoretical results.



\noindent\underline{\hspace*{6.3in}}

\vfil\eject
}
\end{titlepage}

\sect{Introduction}

The Hestenes-Powell augmented Lagrangian has been attracting extensive attentions from optimization researchers.
Many augmented Lagrangian methods of multipliers, which minimize an augmented Lagrangian approximately and circularly with update of multipliers, have been proposed for various optimization problems, for example, some recent references can be found in \cite{AndBMS08a, BirMar14, CDZ15, HFD16, KS19, LiuLM19}.
We consider the augmented Lagrangian method of multipliers (ALM for short) for solving the nonlinear optimization with only general inequality constraints
\bea \hbox{minimize}\quad(\min)\dd\dd f(x) \label{prob1-1}\\
\hbox{subject to}\quad(\st)\dd\dd c(x)\ge 0, \label{prob1-2} \eea where $x\in\Re^n$, $c(x)=(c_1(x), \ldots, c_m(x))\in\Re^m$, $f: \Re^n\to\Re$ and all $c_i: \Re^n\to\Re (i=1,\ldots,m)$ are supposed to be twice continuously differentiable real-valued functions defined on $\Re^n$ and at least one of these functions is a nonlinear and possibly nonconvex function.
For simplicity, we do not incorporate equality constraints into the problem.
Our method can easily be extended to solve the optimization with general equality and inequality constraints (see the later section on numerical experiments for details). If $f$ is linear and all $c_i (i=1,\ldots,m)$ are affine functions, problem \reff{prob1-1}--\reff{prob1-2} is a dual form of the standard linear programming problem (for example, see \cite{NocWri99,SunYua06}).

The ALM was initially proposed by Hestenes \cite{hesten} and Powell \cite{powell} for solving optimization problems with only equality constraints. It was generalized by Rockafellar \cite{rockaf1} to the optimization problems with inequality constraints. The global and local analysis on ALM has been done by many researchers, for example, \cite{AndBMS08, bertse76, bertse82, ConGoT91, powell, rockaf1, rockaf2, rockaf3, a11}.
Since its birth, the ALM has been playing a very important role in the development of effective numerical methods for convex and nonconvex optimization problems. As its extensions and applications, the alternating direction methods of multipliers (ADMM) originated from the ALM has attracted great attentions.

The state-of-the-art solver LANCELOT (see \cite{ConGoT92}) is one of well known examples of ALM for optimization with inequality constraints. By introducing slack variables $z_i (i=1,\ldots,m)$ in the inequality constraints of \reff{prob1-2} and using the augmented Lagrangian function on equality constraints, Conn, Gould and Toint \cite{ConGoT92} solves a sequence of relaxed subproblems of the form \bea \min_{x,z}\dd\dd L_A(x,z,s;\rho)\equiv f(x)-s^T(c(x)-z)+\hf\rho\|c(x)-z\|^2 \label{cgtp}\\
\hbox{s.t.}\dd\dd z\ge 0, \label{cgtpa}\eea
where $z=(z_i)\in\Re^m$, $s\in\Re^m$ is an estimate of the multiplier vector, $\rho>0$ is a penalty parameter. Both $s$ and $\rho$ are held fixed during the solution of each subproblem and are updated adaptively in virtue of the convergence and feasibility of the approximate solution of the subproblem. Noting that problem \reff{cgtp}--\reff{cgtpa} is a nonlinear program with nonnegative constraints, the LANCELOT uses a projected gradient method to solve this problem (see \cite{ConGoT88}). Thanks to the convexity of the objective function with respect to $z$, an alternative to handling the nonnegative constraints in \reff{cgtp} is to eliminate $z$ from $L_A(x,z,s;\rho)$ by using the minimizer $z=\max\{c(x)-{s}/{\rho},0\}$, which results in an unconstrained optimization subproblem \bea
\min_x\ f(x)+\sum_{i=1}^m\phi(c_i(x),s_i;\rho), \label{cgtp1}\eea
where $$\phi(c_i(x),s_i;\rho)=\left\{\ba{cl}
-s_ic_i(x)+\hf\rho c_i(x)^2, & c_i(x)\le s_i/\rho; \\
-\hf s_i^2/\rho, & c_i(x)>s_i/\rho. \ea\right.$$
Unfortunately, the function $\phi$ in \reff{cgtp1} is in general discontinuous in the second derivative with respect to $x$, see \cite{NocWri99} for detailed discussion.

\subsection{Our contributions.}
In general, a generalization of the ALM for optimization with inequality constraints often depends on some kinds of projections with respect to either primal or dual variables (such as it may request either the primal or dual variables to be nonnegative), which may bring about that the augmented Lagrangian of the ALM for optimization with inequality constraints is generally not twice differentiable. Moreover, the global and local convergence on the ALM is generally established on the condition that the penalty parameter $\rho$ is sufficiently large.

Our main contributions in this paper are as follows. Firstly, we introduce a new twice differentiable augmented Lagrangian for optimization with inequality constraints. Our new augmented Lagrangian is a generalization of the Hestenes-Powell augmented Lagrangian and a combination of the augmented Lagrangian and the interior-point technique. A strict local minimizer of the original problem is an approximate strict local solution of our augmented Lagrangian. Secondly, we present a novel augmented Lagrangian method of multipliers for optimization with inequality constraints. Our method shares a similar algorithmic framework with the existing ALMs for optimization with inequality constraints, but it can use the second derivatives and does not depend on projections on the set of inequality constraints. The dual iterates, penalty and smoothing parameters are updated adaptively.
Thirdly, the global and local convergence of our ALM is established without assuming any constraint qualification and without requiring that the penalty parameter $\rho$ is sufficiently large. It is proved that the proposed method is of strong global convergence. In particular, we have either the results that $\rho_k\to\infty$, $\mu_k$ is far from zero, and any cluster point of the generated sequence is an infeasible stationary point of the original problem, or the results that $\rho_k$ is bounded above, $\mu_k\to 0$, and any cluster point of the generated sequence is a KKT point. If $\rho_k\to\infty$ and $\mu_k\to 0$, then there is a cluster point of the generated sequence which is a singular stationary point of the original problem.
Locally, it is proved that our method is of superlinear/quadratic convergence when it converges to an infeasible stationary point, which shows that our method is capable of rapidly detecting the infeasibility. Under suitable conditions, our algorithm is able to converge in at least linear rate to the KKT point. This is a result consistent with the ALM for optimization with equality constraints. The preliminary numerical experiments on some small benchmark problems demonstrate our theoretical results.

\subsection{Some related works.}
Combining the interior-point technique with the augmented Lagrangian has been a very useful approach in developing effective methods for inequality constrained nonlinear programs in the literature. For example, Goldfarb et al. \cite{GolPSY99} have considered the combination of the augmented Lagrangian and the logarithmic-barrier method. Their method is based on a so-called MBAL function which treats inequality constraints with a modified barrier term and equalities
with an augmented Lagrangian term. The MBAL method alternatively minimizes the MBAL function in the primal
space and updates the Lagrange multipliers.
Very recently, Gill et al. \cite{GilKuR20} presents a primal-dual shifted penalty-barrier method for solving problem \reff{prob1-1}--\reff{prob1-2} by combining the interior-point technique with the augmented Lagrangian. The method is proposed based on minimizing a shifted primal-dual penalty-barrier function. It is shown that a limit point of the sequence of iterates
may always be found that is either an infeasible stationary point or a complementary approximate KKT point which is a KKT point under a regularity
condition that is the weakest constraint qualification associated with sequential
optimality conditions, and the method can be equivalent to a shifted variant of the primal-dual
path-following method in the neighborhood of a solution.

The focus on the strong global convergence of algorithms for nonlinear programs is based on the fact that a minimizer of constrained optimization problems may not be a KKT point but is a singular stationary point (see the test problem (TP5) in section 6). In addition, when the solved optimization problem is infeasible, it is meaningful to know what is the minimal $\ell_1$ or $\ell_2$ norm measurement of constraint violations, which has been one of hot topics in optimization, the interested readers can refer to \cite{BurCuW14,BurHan89,ByrCuN10,ByrMaN01,CheGol06,GilKuR20,LPS,LiuYua08,NocOzW12,yuan95}. In particular, some of these works, such as \cite{BurCuW14,ByrCuN10,DLS17}, have shown that their methods have rapid convergence to the infeasible stationary point.

\subsection{Organization and notations.}
Our paper is organized as follows. In section 2, we propose a twice differentiable augmented Lagrangian for optimization with general inequality constraints. Our new augmented Lagrangian can be considered as a generalization of the Hestenes-Powell augmented Lagrangian for inequality constrained optimization. In order to find a matching dual estimate, a particular minimax problem is introduced. Based on this minimax problem, we present our augmented Lagrangian method of multipliers for nonlinear programs with inequality constraints in section 3. We prove the global and local convergence results of our method for nonlinear programs in sections 4 and 5, respectively. The discussions on extension to solving the nonlinear programming with inequality and equality constraints are presented in section 6. Some numerical results are also reported in this section. We conclude our paper in the last section.

Throughout the paper, we use standard notations from the literature. A letter with
subscript $k$ is related to the $k$th
iteration, the subscript $i$ indicates the $i$th component of a vector, and the subscript $kj$ is the $j$th
iteration in solving the $k$th subproblem. All vectors are column vectors, and $u=(x,s)$ means $u=[x^T,\hspace{2pt}s^T]^T$, $s\in\Re_{++}^m$ means $s_i>0$ for $i=1,\ldots,m$. The expression
$\theta_k={O}(t_k)$ means that there exists a scalar $M$
independent of $k$ such that $|\theta_k|\le M|t_k|$ for all $k$ large enough, and
$\theta_k={o}(t_k)$ indicates that $|\theta_k|\le\epsilon_k|t_k|$ for all $k$ large enough with $\lim_{k\to
0}\epsilon_k=0$. If it is not specified, $I$ is the identity matrix and ${\cal I}$ is an index set, $\|\cdot\|$ is the Euclidean norm. Some unspecified notations may be
identified from the context.

\sect{A twice differentiable augmented Lagrangian}

By introducing a slack vector $z\in\Re^m$, the problem \reff{prob1-1}--\reff{prob1-2} can be reformulated as the one with general equality and nonnegative constraints \bea
\min\dd\dd f(x) \label{prob2-1}\\
\st\dd\dd c(x)-z=0,\label{prob2-2}\\
\dd\dd z\ge 0. \label{prob2-3}\eea
It is well known that the original problem \reff{prob1-1}--\reff{prob1-2} and its reformulation \reff{prob2-1}--\reff{prob2-3} are equivalent in the sense that both problems have the same feasible solutions and minimizers, for example, see \cite{BS0,BS,CheGol06,LiuSun01,LiuYua07,NocOzW12,WacBie00}.

Primal-dual interior-point approach has been demonstrated to be very efficient in solving linear and nonlinear constrained optimization problems.
The generic primal-dual interior-point approach for problem \reff{prob1-1}--\reff{prob1-2} generates the interior-point iterates by an inner algorithm approximately solving the unconstrained logarithmic-barrier subproblem \bea
\min_{x}\quad f(x)-\mu\sum_{i=1}^m\ln c_i(x) \label{lgbp}\eea
or its corresponding parametric differential system, where $x$ is requested to be always strictly feasible, $\mu>0$ is a barrier parameter which is held fixed when solving the subproblem \reff{lgbp} or its parametric system.
It is commonly believed that finding a feasible point of the optimization problem with general inequality constraints is almost as difficult as solving the problem.
In order to avoid the difficulty in finding a strict feasible point, we often replace the subproblem \reff{lgbp} by the following equality constrained logarithmic-barrier subproblem for the reformulation \reff{prob2-1}--\reff{prob2-3},
\bea \min_{x,z}\dd\dd f(x)-\mu\sum_{i=1}^m\ln z_i \label{lgbp1}\\
\hbox{s.t.}\dd\dd c(x)-z=0, \label{lgbp1a}\eea
where $z>0$ should be kept for the whole iterative process. However, problem \reff{prob2-1}--\reff{prob2-3} and its associated subproblem \reff{lgbp1}--\reff{lgbp1a} have more variables than the subproblem \reff{lgbp}, and, even if the original problem \reff{prob1-1}--\reff{prob1-2} is convex, problem \reff{lgbp1}--\reff{lgbp1a} can be nonconvex when some of $c_i(x) (i=1,\ldots,m)$ are nonlinear.

Generally, the primal-dual interior-point methods for problem \reff{prob1-1}--\reff{prob1-2} solve the subproblem \reff{lgbp1}--\reff{lgbp1a} or its associated KKT system approximately by an inner algorithm, and the slack variables are requested to be uniformly positive during the iterative process. Recent works on interior-point relaxation methods \cite{DLS17, LiuDai18, LiuDaH20} proposed a new approach which does not solve the logarithmic-barrier problem \reff{lgbp1}--\reff{lgbp1a} directly. Suppose that we have a binary function $t=w(a,b)$ which satisfies $t>0$ for some given $b\in\Re$ and all $a\in\Re$, and there exists $a>0$ such that $a=w(a,b)$ (some parameters may be included in $w$). Such an example is \bea w(a,b;\mu,\rho)=\frac{1}{2\rho}(\sqrt{(b-\rho a)^2+4\rho\mu}-(b-\rho a)), \label{20210523a}\eea
where $\mu>0$ and $\rho>0$ are two given parameters. The function \reff{20210523a} is derived from a closed-form solution on slack variables $z_i\ (i=1,\ldots,m)$ of the stationary conditions of augmented Lagrangian of the logarithmic-barrier problem \reff{lgbp1}--\reff{lgbp1a}, \bea
z_i=\frac{1}{2\rho}(\sqrt{(s_i-\rho c_i(x))^2+4\rho\mu}-(s_i-\rho c_i(x))), \label{20210604a}\eea
where $s=(s_i)\in\Re^m$ is an estimate of the vector of multipliers associated with the equality constraints in \reff{lgbp1a}, $\mu>0$ is the barrier parameter and $\rho>0$ is the penalty parameter. Since the logarithmic-barrier problem \reff{lgbp1}--\reff{lgbp1a} is an equality constrained optimization, the Hestenes-Powell augmented Lagrangian is available and has the form \bea
L_B(x,z,s;\mu,\rho)=f(x)-\mu\sum_{i=1}^m\ln z_i-s^T(c(x)-z)+\hf\rho\|c(x)-z\|^2. \eea
Then the expression of $z_i$ in \reff{20210604a} for $i=1,\ldots,m$ follows from $\na_{(x,z,s)}L_B(x,z,s;\mu,\rho)=0$.
For more details on the derivation of the formulae of $z_i$, the interested readers can refer to \cite{DLS17,LiuDai18}.
Using the binary function $w: \Re\times\Re\to\Re$, instead of solving \reff{lgbp1}--\reff{lgbp1a}, we consider the following modified logarithmic-barrier problem \bea
\min_{x\in\Re^n}\dd\dd f(x)-\mu\sum_{i=1}^m\ln w(c_i(x), s_i; \mu, \rho) \label{lgbp2a}\\
\st\dd\dd c_i(x)-w(c_i(x), s_i; \mu, \rho)=0,\ i=1,\ldots,m. \label{lgbp2b}
\eea
\ble\label{lem0604a} Suppose that $\mu>0$ and $\rho>0$, $s\in\Re^m$ is given. If $(x^*,z^*)$ is a local solution of the logarithmic-barrier subproblem \reff{lgbp1}--\reff{lgbp1a} and $c_i(x^*)=w(c_i(x^*),s_i;\mu,\rho)$ for $i=1,\ldots,m$,
then $x^*$ is a local solution of the modified problem \reff{lgbp2a}--\reff{lgbp2b}.
In addition, if $((x^*,z^*),s^*)$ is a KKT pair of the logarithmic-barrier subproblem, and the binary function $w$ is defined by \reff{20210523a}, then $x^*$ is a KKT point of problem \reff{lgbp2a}--\reff{lgbp2b} with $s=s^*$, and $s_i^*$ is the associated Lagrange multiplier with the constraint $c_i(x)-w(c_i(x),s_i^*;\mu,\rho)=0$ for $i=1,\ldots,m$. \ele\prf
If $x^*$ is a local solution of \reff{lgbp1}--\reff{lgbp1a}, then $z^*=c(x^*)>0$ and \bea
 f(x^*)-\mu\sum_{i=1}^m\ln z_i^*=f(x^*)-\mu\sum_{i=1}^m\ln c_i(x^*)\le f(x)-\mu\sum_{i=1}^m\ln c_i(x) \nn\eea for all $x$ such that $c(x)>0$. Thus,
\bea f(x^*)-\mu\sum_{i=1}^m\ln w(c_i(x^*),s_i;\mu,\rho)\le f(x)-\mu\sum_{i=1}^m\ln w(c_i(x),s_i;\mu,\rho) \nn\eea
for all $x$ such that, for $i=1,\ldots,m$, $c_i(x)=w(c_i(x),s_i;\mu,\rho)$. Thus, $x^*$ is a local solution of the modified problem \reff{lgbp2a}--\reff{lgbp2b}.

For given parameters $\mu>0$ and $\rho>0$, let
$$L(x,s,\lambda)=f(x)-\mu\sum_{i=1}^m\ln w(c_i(x), s_i; \mu, \rho)-\sum_{i=1}^m\lambda_i(c_i(x)-w(c_i(x), s_i; \mu, \rho))$$
be the Lagrange function of problem \reff{lgbp2a}--\reff{lgbp2b}.
Then \bea
\na_x L(x,s,\lambda)=\dd\dd\na f(x)-\sum_{i=1}^m[(\mu/w(c_i(x), s_i; \mu, \rho))w'(c_i(x), s_i; \mu, \rho) \nn\\
\dd\dd+\lambda_i(1-w'(c_i(x), s_i; \mu, \rho))]\na c_i(x). \label{20210529a}\eea
If $((x^*,z^*),s^*)$ is a KKT pair of the logarithmic-barrier subproblem \reff{lgbp1}--\reff{lgbp1a}, then
$\na f(x^*)-\sum_{i=1}^m s_i^*\na c_i(x^*)=0$, $s_i^*=\mu/c_i(x^*)$ and $c_i(x^*)-z_i^*=0$ for $i=1,\ldots,m$. Therefore, due to \reff{20210523a},
$c_i(x^*)=w(c_i(x^*),s_i^*;\mu,\rho)$ for $i=1,\ldots,m$, and, by \reff{20210529a}, $\na_x L(x^*,s^*,s^*)=0$.
\eop

Note that the modified problem \reff{lgbp2a}--\reff{lgbp2b} is an optimization problem with equality constraints, we can similarly use the Hestenes-Powell augmented Lagrangian to reformulate it as the unconstrained optimization
\bea \min_{x\in\Re^n}\ F(x,s;\mu,\rho)=f(x)+\sum_{i=1}^m\psi(c_i(x),s_i;\mu,\rho), \label{20210529b}\label{newalm}\eea
where \bea\dd\dd\psi(c_i(x),s_i;\mu,\rho)\nn\\
\dd\dd=-\mu\ln w(c_i(x),s_i;\mu,\rho)-s_i(c_i(x)-w(c_i(x),s_i;\mu,\rho))+\hf\rho (c_i(x)-w(c_i(x),s_i;\mu,\rho))^2. \nn\eea
Inspired by \refl{lem0604a}, we still use the $s$ in the modified problem as the estimate of the Lagrange multiplier vector of the augmented Lagrangian \reff{20210529b}. Thus, it is a generalization of the classic Hestenes-Powell augmented Lagrangian and is a combination of the augmented Lagrangian and the interior-point technique.

Although the subproblem \reff{newalm} is similar to the augmented Lagrangian subproblem \reff{cgtp1} and the logarithmic-barrier subproblem \reff{lgbp} in appearance that all of them are unconstrained optimization and first-order smooth, but it is essentially distinct from the latter two subproblems in the following aspects. \begin{itemize}
\item[(1)] The function $\psi(c_i(x),s_i;\mu,\rho)$ in \reff{newalm} has one more parameter $\mu$ than $\phi(c_i(x),s_i;\rho)$ in \reff{cgtp1}, $\psi$ is always twice continuously differentiable with respect to $x$ provided $c_i$ is twice continuously differentiable and $s_i$ holds fixed, while $\phi$ of \reff{cgtp1} has discontinuous second derivative with respect to $x$.

\item[(2)] The subproblems \reff{newalm} and \reff{lgbp} are convex if the original problem \reff{prob1-1}--\reff{prob1-2} is convex, while the equivalent problem \reff{lgbp1}--\reff{lgbp1a} of subproblem \reff{lgbp} can be nonconvex even though the original problem is convex.

\item[(3)] Unlike subproblems \reff{lgbp} and \reff{lgbp1}--\reff{lgbp1a}, subproblem \reff{newalm} does not require either $x$ to be feasible or $s$ to be nonnegative. Moreover, $\psi(c_i(x),s_i;\mu,\rho)$ is well-defined for every $x\in\Re^n$ and $s\in\Re^m$, while \reff{lgbp} and \reff{lgbp1}--\reff{lgbp1a} request $c(x)>0$ and $z>0$, respectively.
\end{itemize}

For convenience of statement, we follow our previous works \cite{LiuDai18, LiuDaH20} to similarly define
$z:\Re^{n+m}\to\Re^m$ and $y:\Re^{n+m}\to\Re^m$ to be functions on $(x,s)$ by components \bea
\dd\dd z_i(x,s;\mu,\rho)\equiv\frac{1}{2\rho}(\sqrt{(s_i-\rho c_i(x))^2+4\rho\mu}-(s_i-\rho c_i(x))), \label{zydf1}\\
\dd\dd y_i(x,s;\mu,\rho)\equiv\frac{1}{2\rho}(\sqrt{(s_i-\rho c_i(x))^2+4\rho\mu}+(s_i-\rho c_i(x))), \label{zydf}\eea
where $i=1,\ldots,m$, $x\in\Re^n$ and $s=(s_i)\in\Re^m$ are variables, $\mu>0$ and $\rho>0$ are given parameters. In this writing, \bea
\dd\dd w(c_i(x), s_i;\mu,\rho)=z_i(x,s;\mu,\rho), \\
\dd\dd\psi(c_i(x),s_i;\mu,\rho)=h_i(x,s;\mu,\rho)\equiv -\mu\ln z_i(x,s;\mu,\rho)+\hf\rho |y_i(x,s;\mu,\rho)|^2-\frac{1}{2\rho}s_i^2, \label{Gforma}\eea
where $h_i:\Re^{n+m}\to\Re$ is a real-valued function.
Throughout the paper, we write $z_i(x,s;\mu,\rho)$ and $y_i(x,s;\mu,\rho)$ as $z_i$ and $y_i$ for simplicity when they are not confusing.

The preliminary results are similar to that in \cite{LiuDaH20} and can be proved similarly.
\ble\label{lemzp} For given $\mu>0$ and $\rho>0$, $z_i$ and $y_i$
are defined by \reff{zydf1} and \reff{zydf}. Then \\
(1) $z_i>0$, $y_i>0$, $\rho(c_i(x)-z_i)=s_i-\rho y_i$, and $\rho z_iy_i=\mu$; \\[5pt]
(2) $c_i(x)>0,\ s_i>0,\ c_i(x)s_i=\mu$ if and only if $c_i(x)-z_i=0$; \\[5pt]
(3) $\rho(z_i+y_i)=\sqrt{(s_i-\rho c_i(x))^2+4\rho\mu}$; \\[5pt]
(4) $\rho(z_i-c_i(x))(y_i+c_i(x))=\mu-c_i(x)s_i$; \\[5pt]
(5) $z_i$ and $y_i$ are differentiable, respectively, with respect to $x$ and $s$, and \bea
&&\na_xz_i=\frac{z_i}{z_i+y_i}\na c_i(x), \quad \na_xy_i=-\frac{y_i}{z_i+y_i}\na c_i(x), \label{20140327a}\\
&&\na_sz_i=-\frac{1}{\rho}\frac{z_i}{z_i+y_i}e_i, \quad
\na_sy_i=\frac{1}{\rho}\frac{y_i}{z_i+y_i}e_i, \label{20140327b} \eea
where $e_i\in\Re^m$ is the $i$-th coordinate vector.  \ele\prf
(1) For any given $\mu$, the condition $\mu>0$ implies that \bea
z_i>\frac{1}{2\rho}(|s_i-\rho c_i(x)|-(s_i-\rho c_i(x)))\ge 0. \nn\eea
Similarly, one has $y_i>0$. Note that $\rho (y_i-z_i)=s_i-\rho c_i(x)$ {and} $4\rho^2 z_iy_i=4\rho\mu$. The result (1) follows immediately.

(2) If $c_i(x)-z_i=0$, then, due to (1), $c_i(x)=z_i>0$ and $s_i-\rho y_i=\rho(c_i(x)-z_i)=0$, which further implies $s_i=\rho y_i>0$. Thus,
$c_i(x)s_i=\rho z_iy_i=\mu$.

If $c_i(x)>0,\ s_i>0,\ c_i(x)s_i=\mu$, we need to prove that $c_i(x)-z_i=0$. Suppose that it is not the case, that is, $c_i(x)\ne z_i$. For example, we may assume $c_i(x)>z_i$. Then, by (1), one has $s_i>\rho y_i$, which results in $c_i(x)s_i>\mu$, a contradiction to the condition $c_i(x)s_i=\mu$. The contradiction shows that the result holds.

(3) The result is straightforward from the definitions \reff{zydf1} and \reff{zydf}.

(4) Since $\rho(z_i-c_i(x))(y_i+c_i(x))=\rho z_iy_i-\rho (y_i-z_i+c_i(x))c_i(x)$, the result follows immediately from (1).

(5) By (1), $\rho (z_i-y_i)=\rho c_i(x)-s_i$. Together with (3), one has \bea
\rho(\na_x z_i-\na_x y_i)=\rho\na c_i(x), \ \rho(\na_x z_i+\na_x y_i)=\frac{\rho c_i(x)-s_i}{z_i+y_i}\na c_i(x). \nn\eea
Thus, by doing summation and subtraction, respectively, on both sides of the preceding equations, we have
\bea \dd\dd 2\rho\na_x z_i=(\rho+\frac{\rho c_i(x)-s_i}{z_i+y_i})\na c_i(x)=\frac{2\rho z_i}{z_i+y_i}\na c_i(x), \nn\\[5pt]
\dd\dd -2\rho\na_x y_i=(\rho-\frac{\rho c_i(x)-s_i}{z_i+y_i})\na c_i(x)=\frac{2\rho y_i}{z_i+y_i}\na c_i(x), \nn \eea
where the last equalities in the preceding two equations are obtained from the fact that $\rho (y_i+c_i(x))=\rho z_i+s_i$. Therefore, \reff{20140327a} follows immediately. The results in \reff{20140327b} can be derived in the same way by differentiating with respect to $s$.
\eop

The preceding lemma shows that both functions $z_i$ and $y_i$ are smooth when $\mu>0$.
The following results show that the augmented Lagrangian $F(x,s;\mu,\rho)$ in \reff{20210529b} plays a key role like a penalty function, in which $\rho$ is the penalty parameter and every $h_i$ is a smooth penalty term promoting the iterate to become strictly feasible to the $i$-th inequality constraint of the original problem.
\ble\label{fgpara}\label{yp}
Given $\mu>0$ and $\rho>0$. Let $z=z(x,s;\mu,\rho)$ and $y=y(x,s;\mu,\rho)$ be defined by \reff{zydf1} and \reff{zydf}, $h_i$ be a function given by \reff{Gforma}. Then \\[5pt]
(1) The function $h_i$ is differentiable with respect to $\rho$, and \bea
\frac{\partial h_i(x,s;\mu,\rho)}{\partial\rho}=\hf(c_i(x)-z_i)^2. \nn\eea
That is, if $c_i(x)-z_i\ne 0$, then $h_i$ is a monotonically increasing function with respect to $\rho$. \\[5pt]
(2) There holds \bea
\frac{\partial (c_i(x)-z_i)^2}{\partial \rho}=-\frac{2}{\rho}\frac{z_i}{z_i+y_i}(c_i(x)-z_i)^2. \label{20190622a}\eea
It shows that $(c_i(x)-z_i)^2$ will be reduced as $\rho$ is increased.
\ele\prf (1) It is known that $h_i$ is differentiable with respect to $\rho$, and \bea
\frac{\partial h_i(x,s;\mu,\rho)}{\partial\rho}=-\frac{\mu}{z_i}\frac{\partial z_i}{\partial\rho}+\hf y_i^2+\rho y_i\frac{\partial y_i}{\partial\rho}+\frac{1}{2\rho^2}s_i^2. \nn\eea
Since $\frac{\partial (\rho z_i)}{\partial\rho}=\frac{z_i(y_i+c_i(x))}{z_i+y_i}$, one has $\frac{\partial z_i}{\partial\rho}=\frac{z_i(c_i(x)-z_i)}{\rho(z_i+y_i)}$.  Similarly, $\frac{\partial y_i}{\partial\rho}=-\frac{y_i(c_i(x)+y_i)}{\rho(z_i+y_i)}$. Then \bea
\frac{\partial h_i(x,s;\mu,\rho)}{\partial\rho} 
\dd\dd=-\frac{\mu(c_i(x)-z_i)+\rho y_i^2(c_i(x)+y_i)}{\rho(z_i+y_i)}+\hf y_i^2+\frac{1}{2\rho^2}s_i^2 \nn\\
\dd\dd=-\frac{z_i(c_i(x)-z_i)+y_i(c_i(x)+y_i)}{z_i+y_i}y_i+\hf y_i^2+\frac{1}{2\rho^2}s_i^2 \nn\\
\dd\dd=(-c_i(x)+z_i-y_i)y_i+\hf y_i^2+\frac{1}{2\rho^2}s_i^2 \nn\\
\dd\dd=\hf(\frac{1}{\rho}s_i-y_i)^2. \nn\eea
Thus, the result follows from the equation $\rho (c_i(x)-z_i)=s_i-\rho y_i$.

(2) Due to $\frac{\partial z_i}{\partial\rho}=\frac{z_i(c_i(x)-z_i)}{\rho(z_i+y_i)}$, \bea
\frac{\partial (c_i(x)-z_i)^2}{\partial \rho}=-2(c_i(x)-z_i)\frac{\partial z_i}{\partial\rho}=-\frac{2}{\rho}\frac{z_i}{z_i+y_i}(c_i(x)-z_i)^2, \nn
\eea which completes the proof. \eop

The subsequent results illustrate that our augmented Lagrangian has the same differentiability as $f(x)$ and $c(x)$, and show the close relation between the original inequality constrained optimization \reff{prob1-1}--\reff{prob1-2} and the unconstrained optimization \reff{20210529b}.

\ble\label{bfpropa} Given $\mu>0$ and $\rho>0$. If $f: \Re^n\to\Re$ and $c: \Re^n\to\Re^m$ are twice continuously differentiable on $\Re^n$, then \\
(1) $F$ is twice differentiable with respect to $x$, and
\bea \dd\dd \na_xF(x,s;\mu,\rho)=\na f(x)-\rho\na c(x) y,\nn\\
\dd\dd\na_x^2F(x,s;\mu,\rho)=\left(\na^2 f(x)-\rho\sum_{i=1}^m y_i\na^2 c_i(x)\right)+\rho\sum_{i=1}^m \frac{y_i}{z_i+y_i}\na c_i(x)\na c_i(x)^T, \nn \eea
where $Z=\diag(z)$ and $Y=\diag(y)$; \\
(2) if $f$ and $-c_i\ (i=1,\ldots,m)$ are convex on $\Re^n$, then $F(x,s;\mu,\rho)$ is a convex function with respect to $x$ on $\Re^n$. \ele\prf (1)
Due to \reff{Gforma}, \bea \na_x h_i(x,s;\mu,\rho)=-\frac{\mu}{z_i}\na_x z_i+\rho y_i\na_x y_i=\frac{-\mu-\rho y_i^2}{z_i+y_i}\na c_i(x)=-\rho y_i\na c_i(x). \nn
\eea Thus, $\na_xF(x,s;\mu,\rho)=\na f(x)+\sum_{i=1}^m\na_x h_i(x,s;\mu,\rho)=\na f(x)-\rho\na c(x)y$. Furthermore, by \refl{yp} (1), \bea
\na_x^2 h_i(x,s;\mu,\rho)=-\rho y_i\na^2 c_i(x)+\rho\frac{y_i}{z_i+y_i}\na c_i(x)\na c_i(x)^T. \nn\eea
Therefore, the expression on $\na_x^2F(x,s;\mu,\rho)$ is obtained since \bea
\na_x^2F(x,s;\mu,\rho)=\na^2 f(x)+\sum_{i=1}^m\na_x^2 h_i(x,s;\mu,\rho). \nn \eea

The result (2) is straightforward since the Hessian of $F(x,s;\mu,\rho)$ with respect to $x$ is positive semi-definite. \eop

The results in \refl{bfpropa} show that, if the original problem \reff{prob1-1}--\reff{prob1-2} is convex, then $F(x,s;\mu,\rho)$ is always convex with respect to $x$ for all positive parameters $\mu$ and $\rho$. This is different from problem \reff{lgbp1}--\reff{lgbp1a}, for which the convexity may be destroyed by the introduction of slack variables.

\bth\label{a} Let $x^*$ be a local minimizer of problem \reff{prob1-1}--\reff{prob1-2} at which the linear independence constraint qualification and the second-order sufficient conditions are satisfied with $s=s^*$. Then for $s=s^*$ and $\mu>0$ sufficiently small, there exists a threshold value $\tilde\rho>0$ independent of $\mu$ such that for all $\rho\ge\tilde\rho$, $x^*$ is a $\sqrt{\rho\mu}$-approximate strict local minimizer of the augmented Lagrangian \reff{20210529b} (that is, there a scalar $\delta>0$ such that $\|\na_x F(x^*,s^*;\mu,\rho)\|\le\delta\sqrt{\rho\mu}$).
\eth\prf Under the conditions of the theorem, $x^*$ is a KKT point of problem \reff{prob1-1}--\reff{prob1-2}. Thus, \bea
\na f(x^*)-\na c(x^*)s^*=0,\ c(x^*)\ge 0,\ s^*\ge 0,\ c(x^*)^Ts^*=0. \label{20210603a}\eea
Let $y^*_i=y_i(x^*,s^*;\mu,\rho)$. Note that $c_i(x^*)s_i^*=0$ for $i=1,\ldots,m$. Then \bea
\rho y_i^*=\left\{\ba{ll}
\hf (\sqrt{(s_i^*)^2+4\rho\mu}+s_i^*), & \hbox{if}\ c_i(x^*)=0,\ s_i^*>0; \\[5pt]
\hf (\sqrt{\rho^2 c_i^2(x^*)+4\rho\mu}-\rho c_i(x^*)), & \hbox{if}\ c_i(x^*)>0,\ s_i^*=0;\\[5pt]
\sqrt{\rho\mu}, &\hbox{otherwise.} \ea\right.\nn\eea
Since $\sqrt{(s_i^*)^2+4\rho\mu}\le s_i^*+2\sqrt{\rho\mu}$ and $\sqrt{\rho^2 c_i^2(x^*)+4\rho\mu}\le\rho c_i(x^*)+2\sqrt{\rho\mu}$, one has \bea
s^*\le\rho y^*\le s^*+\sqrt{\rho\mu},\quad \|\rho y^*-s^*\|_{\infty}\le \sqrt{\rho\mu}. \label{20210603b}\eea

If $\na_x F(\hat x,s^*;\mu,\rho)=0$, and $\na^2_{xx} F(\hat x,s^*;\mu,\rho)$ is positive definite, then $\hat x$ is a strict local minimizer of problem \reff{20210529b}. We will prove the result by showing $\|\na_x F(x^*,s^*;\mu,\rho)\|\le\delta\sqrt{\rho\mu}$ for some scalar $\delta$ and $\na^2_{xx} F(x^*,s^*;\mu,\rho)$ is positive definite for all $\rho$ greater than some scalar $\tilde\rho$. By using \refl{bfpropa}, and \reff{20210603a}, \reff{20210603b}, we have \bea
\|\na_x F(x^*,s^*;\mu,\rho)\|=\|\na f(x^*)-\rho\na c(x^*)y^*\|=\|\na c(x^*)(s^*-\rho y^*)\|\le\sqrt{\rho\mu}\|\na c(x^*)\|_1, \nn \eea
which verifies the first part of the result.

Now we prove the second part of the result by showing that $d^T\na^2_{xx} F(x^*,s^*;\mu,\rho)d>0$ for all nonzero $d\in\Re^n$ and all $\rho>0$ sufficiently large. Let $z_i^*=z_i(x^*,s^*;\mu,\rho)$. Then \bea
\frac{y_i^*}{z_i^*+y_i^*}=\left\{\ba{ll}
\hf (1+s_i^*/\sqrt{(s_i^*)^2+4\rho\mu}), & \hbox{if}\ c_i(x^*)=0,\ s_i^*>0; \\[5pt]
\hf (1-\rho c_i(x^*)/\sqrt{\rho^2 c_i^2(x^*)+4\rho\mu}), & \hbox{if}\ c_i(x^*)>0,\ s_i^*=0;\\[5pt]
\hf, &\hbox{otherwise.} \ea\right.\nn\eea
Therefore, by \refl{bfpropa}, \bea
\dd\dd\na^2_{xx} F(x^*,s^*;\mu,\rho) \nn\\
\dd\dd=(\na^2 f(x^*)-\rho\sum_{i=1}^m y_i^*\na^2 c_i(x^*))+\rho\sum_{i=1}^m \frac{y_i^*}{z_i^*+y_i^*}\na c_i(x^*)\na c_i(x^*)^T \nn\\
\dd\dd=(\na^2 f(x^*)-\sum_{i=1}^m s_i^*\na^2 c_i(x^*))+\sum_{i=1}^m (s_i^*-\rho y_i^*)\na^2 c_i(x^*) \nn\\
\dd\dd\quad+\hf\rho(\sum_{i\in I_1}(1+s_i^*/\sqrt{(s_i^*)^2+4\rho\mu})\na c_i(x^*)\na c_i(x^*)^T+\sum_{i\in I_2}\na c_i(x^*)\na c_i(x^*)^T) \nn\\
\dd\dd\quad+\hf\rho\sum_{i\in I_3}(1-\rho c_i(x^*)/\sqrt{\rho^2 c_i^2(x^*)+4\rho\mu})\na c_i(x^*)\na c_i(x^*)^T, \nn\eea
where $I_1=\{i|c_i(x^*)=0,\ s_i^*>0\}$, $I_2=\{i|c_i(x^*)=0,\ s_i^*=0\}$, $I_3=\{i|c_i(x^*)>0,\ s_i^*=0\}$. If $d^T(\na^2 f(x^*)-\sum_{i=1}^m s_i^*\na^2 c_i(x^*))d>0$ for all nonzero $d\in\Re^n$ satisfying $\na c_i(x^*)^Td=0$, $i\in I_1\cup I_2$, then by \reff{20210603b} and the proof of Theorem 17.5 of \cite{NocWri99}, the result follows easily.  \eop

\sect{A novel augmented Lagrangian method of multipliers}

To ensure that $s$ is a good estimate of Lagrange multiplier vector, we maximize the augmented Lagrangian with respect to $s$ in problem \reff{20210529b}, which results in the following unconstrained minimax problem
\bea\min_{x\in\Re^n}\max_{s\in\Re^m} F(x,s;\mu,\rho)=f(x)+\sum_{i=1}^m h_i(x,s;\mu,\rho). \label{spp1} \eea

In problem \reff{spp1}, both $x$ and $s$ are variables. We prove some results on the differentiability of $F(x,s;\mu,\rho)$ with respect to $s$.
\ble\label{bfprop} Given $\mu>0$ and $\rho>0$, $z=z(x,s;\mu,\rho)$ and $y=y(x,s;\mu,\rho)$ are defined by \reff{zydf1} and \reff{zydf}. Then one has the following results.\\
(1) $F$ is twice differentiable with respect to $s$, and
\bea \dd\dd \na_sF(x,s;\mu,\rho)=z-c(x),\nn\\
     \dd\dd \na_s^2F(x,s;\mu,\rho)=-\frac{1}{\rho}(Z+Y)^{-1}Z, \nn \eea
where $Z=\diag(z)$ and $Y=\diag(y)$.\\
(2) $F(x,s;\mu,\rho)$ is a strictly concave function with respect to $s$ on $\Re^m$. \ele\prf (1)
Note that $\na_s h_i(x,s;\mu,\rho)=-\frac{\mu}{z_i}\na_s z_i+\rho y_i\na_s y_i-\frac{1}{\rho}s_i=(y_i-\frac{1}{\rho}s_i)e_i$ and $\na_s^2 h_i(x,s;\mu,\rho)=-\frac{1}{\rho}\frac{z_i}{z_i+y_i}$. The formulae on $\na_sF(x,s;\mu,\rho)$ and $\na_s^2F(x,s;\mu,\rho)$ are derived immediately from the equation $\rho (z-c(x))=\rho y-s$.

The result (2) is straightforward since the negative Hessian of $F(x,s;\mu,\rho)$ with respect to $s$ is positive definite. \eop

Now we are ready to present our main results in this section. The results show that we can obtain a very well approximated KKT solution of the original problem provided $\mu$ is small enough.
\bth\label{mrs} The following results hold. \\
(1) Given $\mu>0$ and $\rho>0$. Let $(x^*,s^*)\in\Re^{n}\times\Re^n$ be a local solution of the minimax problem \reff{spp1}. Then \bea
\dd \na f(x^*)-\na c(x^*)s^*=0, \dd \label{mimkkt1-1}\label{mimkkt1-2}\\
     \dd c(x^*)-z^*=0, \dd \label{mimkkt1-3}\eea
where $z^*=z(x^*,s^*;\mu,\rho)$. \\
(2) If $(x^*,s^*)$ satisfies conditions \reff{mimkkt1-1}--\reff{mimkkt1-3}, and $\mu=0$ and $\rho>0$, then $(x^*,s^*)$ is a KKT pair of the original problem
\reff{prob1-1}--\reff{prob1-2}. \eth\prf (1) If $(x^*,s^*)$ is a local solution of the minimax problem \reff{spp1}, then \bea
\na_x F(x^*,s^*;\mu,\rho)=0, \ \hbox{and}\ \na_s F(x^*,s^*;\mu,\rho)=0, \nn\eea
which, together with \refl{bfpropa} (1), \refl{bfprop} (1) and the fact that $\rho (c(x^*)-z^*)=s^*-\rho y^*$, imply that equations \reff{mimkkt1-1}--\reff{mimkkt1-3} are satisfied.

(2) For $\mu=0$ and $\rho>0$, $z_i^*=\frac{1}{2\rho}(|s_i^*-\rho c_i(x^*)|-(s_i^*-\rho c_i(x^*)))$ for $i=1,\ldots,m$. Thus, due to equation \reff{mimkkt1-3},
if $s_i^*-\rho c_i(x^*)\ge 0$, then $z_i^*=c_i(x^*)=0$, which further implies $s_i^*\ge 0$; otherwise, $s_i^*-\rho c_i(x^*)<0$, $z_i^*=c_i(x^*)-\frac{s_i^*}{\rho}$, hence $s_i^*=0$ and $c_i(x^*)>0$. In summary, if $(x^*,s^*)$ satisfies the equation \reff{mimkkt1-3}, and $\mu=0$ and $\rho>0$, then, for every $i=1,\ldots,m$, $c_i(x^*)\ge 0$, $s_i^*\ge 0$, and $c_i(x^*)s_i^*=0$. Combining with condition \reff{mimkkt1-1}, $(x^*,s^*)$ is precisely a KKT pair of the original problem \reff{prob1-1}--\reff{prob1-2}. \eop

The next result shows the relation between the subproblem \reff{spp1} and the logarithmic barrier subproblem \reff{lgbp}. It is because we do not request any $c_i(x)$ to be positive before termination of the proposed method (in other word, our method is admitted to being asymptotically strictly feasible as $\mu$ is decreasing), our method exhibits a robust and distinguished behavior in solving nonlinear programs with inequality constraints.
\bth
Given $\rho>0$. If $\mu>0$ and $c_i(x)>0$ for all $i=1,\ldots,m$, then problem \reff{spp1} is reduced to the logarithmic-barrier subproblem \reff{lgbp}.
\eth\prf If $c_i(x)>0$, then $s_i^*=\mu/c_i(x)$ maximizes $h_i(x,s;\mu,\rho)$ since $s_i^*$ is the unique solution of equation $\na_s h_i(x,s;\mu,\rho)=0$, i.e., $z_i^*(x,s^*;\mu,\rho)=c_i(x)$. For every $i=1,\ldots,m$, by substituting $s_i^*$ for $s_i$ in problem \reff{spp1}, we have the subproblem \reff{lgbp} immediately. \eop

We describe our algorithm for problem \reff{prob1-1}--\reff{prob1-2} in this section. In the algorithm, the minimax problem \reff{spp1} is solved alternately. Correspondingly, the parameters $\rho$ and $\mu$ are updated adaptively. Following this approach, our method shares a similar framework to the existing ALM for optimization with inequality and equality constraints. Now we describe our ALM  for problem \reff{prob1-1}--\reff{prob1-2}.

\begin{algorithm}
\caption{A novel augmented Lagrangian method of multipliers for problem \reff{prob1-1}--\reff{prob1-2}}
\label{alg1}
{\small \alglist
\item[Given] $(x_0,s_0)\in\Re^{n}\times\Re_{++}^m$, $H_0\in\Re^{n\times n}$, $\mu_0\in (0,1)$, $\rho_0\in [1,\infty)$, $\epsilon\in (0,\mu_0)$.
    Compute residuals of the KKT conditions
    $$E_{01}=\frac{1}{\rho_0}\|\na f(x_{0})-\na c(x_0)s_0\|_{\infty},\ E_{02}=\frac{1}{\rho_0}\|s_0\circ c(x_0)\|_{\infty},\ E_{03}=\|\max\{0, -c(x_0)\}\|_{\infty},$$
    and the residual of infeasible stationarity $$E_{04}=\|\na c(x_0)\max\{0, -c(x_0)\}\|_{\infty}.$$
    Set $k:=0$.

\item[While] either $\max\{E_{k1},E_{k2},E_{k3}\}<\epsilon$ or both $E_{k3}>\epsilon$ and $E_{k4}<\epsilon$, stop the algorithm.

\item[Step] 1 ({\bf Derive the estimates of primal variables}). For given $s_k\in\Re^m$ and parameters $\mu_k$ and $\rho_k$, starting from $x_k$, solve the smooth unconstrained optimization subproblem \bea {\min}_x\ \frac{1}{\rho_k}F(x,s_k;\mu_k,\rho_{k}) \label{subpro}\eea to get an approximate solution $x_{k+1}$ such that  \bea
\|\na_x F(x_{k+1},s_k;\mu_k,\rho_{k})\|_{\infty}(=\|\na f(x_{k+1})-\rho_k\na c(x_{k+1})\hat y_{k+1}\|_{\infty})\le 0.95 \rho_k\mu_k, \label{subcon}\eea
where $\hat y_{k+1}=y(x_{k+1},s_k;\mu_k,\rho_{k})$. 

\item[Step] 2 ({\bf Obtain the estimates of dual variables}). Evaluate $\hat z_{k+1}=z(x_{k+1},s_{k};\mu_k,\rho_k)$. Set \bea
s_{k+1}=\rho_k\hat y_{k+1}(=s_k+\rho_{k}(\hat z_{k+1}-c(x_{k+1})). \label{supdate}\eea

\item[Step] 3 ({\bf Update the parameters}). Evaluate $\tilde z_{k+1}=z(x_{k+1},s_{k+1};\mu_k,\rho_k)$. Set $\tilde E_{k+1}=\|\tilde z_{k+1}-c(x_{k+1})\|_{\infty}.$ If $\tilde E_{k+1}>0.95\mu_k$, then select $\rho_{k+1}\ge 2\rho_k$, and set $\mu_{k+1}=\mu_k$ and $s_{k+1}=s_k$;

    otherwise, select $\mu_{k+1}\le 0.1\mu_k$ and set $\rho_{k+1}=\max\{\rho_k,\|s_{k+1}\|_{\infty}\}$.

\item[Step] 4 ({\bf Compute the residuals}). Evaluate
    $$E_{k+1,1}=\frac{1}{\rho_{k+1}}\|\na f(x_{k+1})-\na c(x_{k+1})s_{k+1}\|_{\infty},\
    E_{k+1,2}=\frac{1}{\rho_{k+1}}\|s_{k+1}\circ c(x_{k+1})\|_{\infty},$$
    $$E_{k+1,3}=\|\max\{0, -c(x_{k+1})\}\|_{\infty},\ E_{k+1,4}=\|\na c(x_{k+1})\max\{0, -c(x_{k+1})\}\|_{\infty}.$$
Set $k:=k+1$.

\item[End] (while)

\eli}

\end{algorithm}

When functions $f$ and $c$ are twice continuously differentiable, it is known from \refl{bfpropa} that $F(x,s_k;\mu_k,\rho_{k})$ is twice continuously differentiable with respect to $x$, thus problem \reff{subpro} in \refal{alg1} can be solved by all efficient algorithms for smooth unconstrained optimization problems in the literature. Moreover, by \refl{bfpropa} (2), if $f$ is a convex function and $c_i\ (i=1,\ldots,m)$ are concave functions, then problem \reff{subpro} is still a convex problem.

\refal{alg1} can easily be extended to the optimization problems with general inequality and equality constraints by replacing the smooth unconstrained optimization subproblem \reff{subpro} with an equality constrained subproblem. This is reverse to the sequential quadratic programming approach, which is extended from that for equality constrained optimization to the optimization problems with inequality constraints.

Let $\hat E_{k+1}=\|\hat z_{k+1}-c(x_{k+1})\|_{\infty}$. The following result shows that, if $\tilde E_{k+1}>0.95\mu_k$ in Step 3 of \refal{alg1}, then $\hat E_{k+1}>0.95\mu_k$ since $\tilde E_{k+1}\le\hat E_{k+1}$. Thus, at every iterate, we have either the case (1) in which $s_{k+1}=s_k$, $\rho_{k+1}\ge2\rho_k$, $\mu_{k+1}=\mu_k$, $\hat E_{k+1}>0.95\mu_k$, or the case (2) where $s_{k+1}$ is derived by \reff{supdate}, $\rho_{k+1}\ge\|s_{k+1}\|_{\infty}$ and either $\rho_{k+1}=\rho_k$ or $\rho_{k+1}>\rho_k$, $\mu_{k+1}\le0.1\mu_k$, and $\tilde E_{k+1}\le 0.95\mu_k$.
\ble
Given $\rho>0$ and $\mu>0$. For any $x\in\Re^n$ and $s\in\Re^m$, if one has $s(x)=s+\rho(z(x,s;\mu,\rho)-c(x))$, then \bea
\|z(x,s(x);\mu,\rho)-c(x)\|\le\|z(x,s;\mu,\rho)-c(x)\|. \eea
In particular, if $|z_i(x,s;\mu,\rho)-c_i(x)|\ne 0$ for some $i=1,\ldots,m$, then \bea
|z_i(x,s(x);\mu,\rho)-c_i(x)|<|z_i(x,s;\mu,\rho)-c_i(x)|. \eea
\ele\prf By \refl{lemzp} (1), $s(x)=\rho y(x,s;\mu,\rho)$, and \bea \dd\dd\|z(x,s(x);\mu,\rho)-c(x)\|=\frac{1}{\rho}\|\rho y(x,s(x);\mu,\rho)-s(x)\|=\|y(x,s(x);\mu,\rho)-y(x,s;\mu,\rho)\|. \nn\eea
Moreover, for $i=1,\ldots,m$, \bea \dd\dd 2\rho(z_i(x,s(x);\mu,\rho)-c_i(x)) \nn\\
\dd\dd=2\rho(y_i(x,s(x);\mu,\rho)-s(x)_i) \nn\\
\dd\dd=2\rho(y_i(x,s(x);\mu,\rho)-y_i(x,s;\mu,\rho)) \nn\\[5pt]
\dd\dd=\sqrt{(s(x)_i-\rho c_i(x))^2+4\rho\mu}-\sqrt{(s_i-\rho c_i(x))^2+4\rho\mu}
+(s(x)_i-s_i) \label{20201031a}\\[5pt]
\dd\dd=(s(x)_{i}-s_{i})\left\{1+\frac{(s(x)_{i}-\rho c_i(x))+(s_{i}-\rho c_i(x))}{\sqrt{(s(x)_{i}-\rho c_i(x))^2+4\rho\mu}
+\sqrt{(s_{i}-\rho c_i(x))^2+4\rho\mu}}\right\}, \nn\eea
and $s(x)_{i}-s_{i}=\rho(z_i(x,s;\mu,\rho)-c_i(x))$. Since the absolute value of the second term inside the brackets of the last equation of \reff{20201031a} is less than $1$, the desired result follows immediately.  \eop

\refal{alg1} suggests that $\{\rho_k\}$ is a monotonically non-decreasing sequence of scalars and $\{\mu_k\}$ is a monotonically non-increasing sequence of scalars. By the update rule of $\rho_k$ (see Step 3), if \refal{alg1} does not terminate finitely, then one of the following three cases will happen: one case is that $\rho_k\to\infty$ as $k\to\infty$ and $\mu_k$ keeps to be a constant after a finite number of iterations, the other case is that parameters $1/\rho_k$ and $\mu_k$ reduce alternately in every finite number of iterations and finally $\rho_k\to\infty$ and $\mu_k\to 0$ as $k\to\infty$, and the another case is that $\mu_k\to 0$ as $k\to\infty$ and $\rho_k$ keeps to be a constant after a finite number of iterations.

According to \refal{alg1}, $s_{k+1}\ge 0$ for all $k\ge 0$. Moreover, the sequence $\{\|s_{k+1}\|/\rho_{k+1}\}$ is always bounded.
In addition, if $\mu_k$ is reduced for some $k$, then $E_{k+1,1}\le 0.95\mu_{k}$.

\sect{Global convergence}

For doing global convergence analysis, we set $\epsilon=0$ in \refal{alg1}. In this situation, the algorithm will not terminate in a finite number of iterations. We firstly consider the case that $\rho_k\to\infty$, and $\mu_{k+1}$ and $s_{k+1}$ keeps unchange for all sufficiently large $k$. In this case, $\tilde E_{k+1}>0.95\mu_k$ for all sufficiently large $k$. After that, we consider the cases in that the condition $\tilde E_{k+1}\le 0.95\mu_k$ is always attained in one or more finite iterations, which can be stated as the following: \\
(1) $\tilde E_{k+1}>0.95\mu_{k}$ and $\tilde E_{k'+1}\le 0.95\mu_{k'}\ (k'>k)$ alternately emerge for every finite iterations; \\
(2) $\tilde E_{k+1}\le 0.95\mu_{k}$ for all sufficiently large $k$. \\
Subsequently, we will analyze the convergence regarding the preceding three cases, respectively.

We need the following blanket assumptions for our global convergence analysis.
\bas\label{ass1}\ \\
(1) The functions $f$ and $c_i\ (i=1,\ldots,m)$ are twice continuously
differentiable on $\Re^n$; \\
(2) The iterative sequence $\{x_k\}$ is in an open bounded set of $\Re^n$.
\eas

Since $\{\|s_{k}\|/\rho_k\}$ is bounded, \refa{ass1} implies that all sequences $\{\tilde z_k\}$, $\{\hat z_k\}$, $\{\tilde y_k\}$, and $\{\hat y_k\}$ are bounded. Our first convergence result focuses on the case that $\tilde E_{k+1}>0.95\mu_{k}$ for all sufficiently large $k$. In this case, the algorithm will converge to an infeasible stationary point of the original problem.
\ble\label{lem4.2} Under \refa{ass1}, if the condition \reff{subcon} is always satisfied,
$\tilde E_{k+1}>0.95\mu_{k}$ for all sufficiently large $k$, then $\rho_k\to\infty$ as $k\to\infty$ and $\mu_k=\mu_{k_0}>0$ for some $k_0>0$, and any cluster point $x^*$ of sequence $\{x_k\}$ is an infeasible point to the problem \reff{prob1-1}--\reff{prob1-2} and satisfies \bea
\na c(x^*)\max\{-c(x^*),0\}=0. \eea
That is, $x^*$ is an infeasible stationary point of the problem \reff{prob1-1}--\reff{prob1-2}.
\ele\prf Without loss of generality, suppose that $\tilde E_{k+1}>0.95\mu_k$ for all $k\ge k_0$. Then, for all $k\ge k_0$, $s_k=s_{k_0}$, $\mu_k=\mu_{k_0}$, $\rho_{k+1}\ge 2\rho_k$.
Thus, $\rho_k\to\infty$, and, for any cluster point $x^*$ of $\{x_k\}$, if $\lim_{k\in{\cal K}, k\to\infty}x_{k}=x^*$, then \bea\lim_{k\in{\cal K}, k\to\infty}\hat z_{k}=\max\{c(x^*),0\},\quad \lim_{k\in{\cal K}, k\to\infty}\hat y_{k}=\max\{-c(x^*),0\}. \eea
Note that the condition $\tilde E_{k+1}>0.95\mu_{k_0}$ for all $k\ge k_0$ implies that $\|\max\{-c(x^*),0\}\|\ne 0$, i.e., $x^*$ is an infeasible point to the problem \reff{prob1-1}--\reff{prob1-2}.
Finally, the desired result follows immediately from dividing $\rho_k$ and taking the limit $k\to\infty$ for $k\in{\cal K}$ on both sides of \reff{subcon}.
\eop

In the latter two cases, the algorithm may converge to either a singular stationary point (in this case $\rho_k\to\infty$) or a KKT point of the original problem.
\ble\label{lem4.4} Under \refa{ass1}, if the condition \reff{subcon} is always satisfied, and there have
$\tilde E_{k+1}>0.95\mu_k$ at some iterates $k$ and $\tilde E_{k'+1}\le 0.95\mu_{k'}$ at some other iterates $k'\ (k'>k, k\to\infty, k'\to\infty)$, then $\rho_k\to\infty$ and $\mu_k\to0$  as $k\to\infty$, and there are a cluster point $x^*$ of sequence $\{x_k\}$ and an associated vector $v^*\in\Re^m$ such that \bea
\na c(x^*)v^*=0,\ v^*\ge 0,\ c(x^*)\ge 0,\ (v^*)^Tc(x^*)=0. \label{210408b}\eea
That is, $x^*$ is a feasible point and a singular stationary point of the problem \reff{prob1-1}--\reff{prob1-2}. \ele\prf
Note that the sequence $\{\rho_k\}$ is monotonically non-decreasing and there is a subsequence of $\{\rho_k\}$ which is strictly increasing. Thus, $\rho_k\to\infty$ as $k\to\infty$. Similarly, one has $\mu_k\to 0$ as $k\to\infty$ since the sequence $\{\mu_k\}$ is monotonically non-increasing and there is a subsequence of $\{\mu_k\}$ which is strictly decreasing.

The conditions imply that there exists an infinite subsequence $\{(x_{k_j+1},s_{k_j+1})\}$ such that $\tilde E_{k_j+1}\le 0.95\mu_{k_j}$ holds for all $k_j>0$ and $\mu_{k_j}\to 0$. Thus, $\lim_{k_j\to\infty}\tilde E_{k_j+1}=0$, which suggests that all cluster points of the subsequence $\{x_{k_j+1}\}$ are feasible points. If $v^*$ is a cluster point of $\{\tilde y_{k_j+1}\}$ and $z^*$ is one of $\{\tilde z_{k_j+1}\}$, then $v^*\ge 0$, $c(x^*)\ge 0$, and $(v^*)^Tc(x^*)=0$ since  \bea z^*-c(x^*)=\lim_{k_j\to\infty}\tilde E_{k_j+1}=0\ \hbox{and}\ (v^*)^Tz^*=\lim_{k_j\to\infty}\tilde y_{k_j+1}^T\tilde z_{k_j+1}=0. \nn\eea

Due to $$\lim_{k_j\to\infty}\|\tilde y_{k_j+1}-s_{k_j+1}/\rho_{k_j+1}\|=\lim_{k_j\to\infty}\tilde E_{k_j+1}=0,$$ the sequence $\{s_{k_j+1}/\rho_{k_j+1}\}$ has the same cluster points as $\{\tilde y_{k_j+1}\}$.
Without loss of generality, suppose that $\lim_{k_j\in{\cal K}, k_j\to\infty} s_{k_j+1}/\rho_{k_j+1}=v^*$.
It follows from \reff{subcon} and \reff{supdate} that one has \bea
\|\na f(x_{k_j+1})-\na c(x_{k_j+1})s_{k_j+1}\|_{\infty}\le 0.95 \rho_{k_j+1}\mu_{k_j+1}. \eea
Dividing $\rho_{k_j+1}$ and taking the limit $k_j\to\infty$ for $k_j\in{\cal K}$ on both sides of the preceding inequality, one has \bea
\|\na c(x^*)v^*\|=0, \nn\eea
which completes our proof. \eop

The following result shows that, under suitable conditions, our algorithm will converge to a KKT pair of the original problem.
\ble\label{lem4.5} Under \refa{ass1}, if the condition \reff{subcon} is always satisfied, $\tilde E_{k+1}\le 0.95\mu_k$ for all sufficiently large $k$, then $\mu_k\to 0$ as $k\to\infty$. In addition, if there is a positive integer $k_0$ such that $\rho_k=\rho_{k_0}$ for all $k\ge k_0$, then every cluster point $x^*$ of sequence $\{x_k\}$ is a KKT point of the original problem \reff{prob1-1}--\reff{prob1-2}. Otherwise, $\rho_k\to\infty$ as $k\to\infty$, and there exists a cluster point of $\{x_k\}$ which is a singular stationary point of the original problem.
\ele\prf If $\tilde E_{k+1}\le 0.95\mu_k$ for all $k\ge k_0$, then, by Step 3 of \refal{alg1}, $\rho_{k+1}\ge\rho_{k_0}$ and $\mu_{k+1}\le 0.1\mu_k$ for all $k\ge k_0$. Thus, $\mu_k\to 0$ as $k\to\infty$. This result together with conditions \reff{subcon} and $\tilde E_{k+1}\le 0.95\mu_k$ suggests \bea\lim_{k\to\infty}\|\tilde z_{k+1}-c(x_{k+1})\|=0,\quad \lim_{k\to\infty}\frac{1}{\rho_{k+1}}\|\na f(x_{k+1})-\na c(x_{k+1})s_{k+1}\|=0. \label{210408a}\eea
If $\rho_{k+1}=\rho_{k_0}$ for all $k\ge k_0$, then $\{s_{k+1}\}$ is bounded and the preceding equations imply that every cluster point $(x^*,s^*)$ of $\{(x_k,s_k)\}$ is a KKT pair of  the original problem \reff{prob1-1}--\reff{prob1-2}. Otherwise,
$\{s_{k+1}\}$ is unbounded, $\rho_{k+1}\to\infty$ as $k\to\infty$, and there exists an infinite subsequence of $\{s_{k+1}/\rho_{k+1}\}$ with $\|s_{k+1}\|_{\infty}/\rho_{k+1}=1$. Then \reff{210408a} implies that there exists a cluster point $x^*$ of $\{x_k\}$ which is a feasible point of the original problem such that \reff{210408b} holds for some $v^*$. \eop

For given $s_k$, $\mu_k,\rho_k$, let us consider the line search methods for problem \reff{subpro}: \bea x_{k,j+1}=x_{k,j}+\alp_{k,j}d_{k,j}, \label{lsch}\eea where $d_{k,j}=-\frac{1}{\rho_k}B_{k,j}^{-1}\na_x F(x_{k,j},s_k;\mu_k,\rho_k)$, for simplicity $B_{k,j}\in\Re^{n\times n}$ is assumed to be positive definite with a uniformly bounded condition number for all $j$, $\alp_{k,j}$ is selected to be maximal in $(0,1]$ such that \bea
F(x_{k,j}+\alp_{k,j}d_{k,j},s_k;\mu_k,\rho_k)\le F(x_{k,j},s_k;\mu_k,\rho_k)+\sigma\alp_{k,j}\na_xF(x_{k,j},s_k;\mu_k,\rho_k)^Td_{k,j}, \nn\eea
where $\sigma\in (0,1)$ is a scalar.
\ble\label{lem4.6} Suppose that functions $f$ and $c_i\ (i=1,\ldots,m)$ are twice differentiable, and $f$ is bounded below on $\Re^n$.
If problem \reff{subpro} is solved by the preceding line search method \reff{lsch}, then the method \reff{lsch}
will be terminated finitely to satisfy the condition \reff{subcon}.
\ele\prf For any given $s$, function $h_i$ denoted by \reff{Gforma} is always bounded below on $\Re^n$. If $f$ is bounded below, then $F(x,s_k;\mu_k,\rho_k)$ is bounded below. Since $\na_x F(x,s_k;\mu_k,\rho_k)$ is always Lipschitz continuous on a bounded open set, by the Theorem 3.2 of \cite{NocWri99}, for any given $\epsilon>0$, there is an iteration $j+1$ such that $\|\na_x F(x_{k,j+1},s_k;\mu_k,\rho_k)\|<\epsilon$.   \eop

In a summary, we have the following global convergence results on \refal{alg1}.
\bth\label{gth} Under \refa{ass1}, if $f$ is bounded below on $\Re^n$, \reff{subcon} is always satisfied, $\epsilon>0$, then our algorithm will terminates finitely at either an approximate KKT point of the original problem, or either an approximate infeasible stationary point or an approximate singular stationary point of original problem.
\eth\prf These results follow from Lemmas \ref{lem4.2}, \ref{lem4.4}, and \ref{lem4.5} immediately. \eop

\sect{Local convergence}

We will analyze the local convergence of our algorithm in this section.
The following blanket assumptions are requested for local convergence analysis.
\bas\label{ass2} \ \\
(1) The sequence $\{x_k\}$ is convergent, i.e., $x_k\to x^*$ as $k\to\infty$; \\
(2) The functions $f$ and $c_i\ (i=1,\ldots,m)$ are twice differentiable on $\Re^n$, and their second derivatives are
Lipschitz continuous at some neighborhood of $x^*$;\\
(3) The gradients $\na c_i(x^*) \ (i\in{\cal I})$ are linearly independent, where ${\cal I}=\{i|c_i(x^*)=0, i=1,\ldots,m\}$.
\eas

Conditions (1)--(2) in \refa{ass2} are commonly used in local convergence analysis for nonlinear programs. Under \refa{ass2} (3), the limit $x^*$ cannot be a singular stationary point. Thus, based on our global convergence analysis, we focus on the local convergence in the following two cases: \\
(1) $x^*$ is an infeasible stationary point, in which case, $\rho_k\to\infty$ as $k\to\infty$, and without loss of generality, we can assume that $s_k=s_0$, $\mu_k=\mu_0>0$ and $\hat E_k>0.95\mu_0$ for all $k>0$; \\
(2) $s_k\to s^*$ and $(x^*,s^*)$ is a KKT pair, in which case, $\mu_k\to 0$ as $k\to\infty$, and without loss of generality, we assume that $\rho_k=\rho_0>0$ and $\tilde E_k\le 0.95\mu_{k-1}$ for all $k>0$.

We will show that our algorithm has the potential of rapidly detecting the possible infeasibility of the problem  \reff{prob1-1}--\reff{prob1-2}. Moreover, under the suitable conditions, our algorithm can be linearly convergent to the KKT point when the original problem is feasible, which is a result similar to that on the ALM for optimization with equality constraints.

\subsection{Rapid convergence to an infeasible stationary point.}
The following assumption is natural on a Newton's method for the nonsmooth equation \bea
\na c(x)\max\{0, -c(x)\}=0, \nn\eea
for example, see \cite{QiSun93}.
\bas\label{ass21} The limit $x^*$ is an infeasible stationary point and matrix
$$B^*=-\sum_{i=1}^m\max\{0,-c_i(x^*)\}\na^2c_i(x^*)+\sum_{i\in\{i|c_i(x^*)<0\}}\na c_i(x^*)\na c_i(x^*)^T,$$ is positive definite. \eas

In virtue of \refa{ass21} and our global convergence results, one has $\rho_k\to\infty$ as $k\to\infty$. Let $x_{kj}$ be the $j$-th iterate generated by the line search methods for minimizing the $k$-th subproblem \reff{subpro}, $B_{kj}$ is the associated Hessian $\frac{1}{\rho_k}\na^2_xF(x_{kj},s_0;\mu_0,\rho_k)$. Then, by the continuities of $\na^2f$, $\na^2c_i$ and $\na c_i (i=1,\ldots,m)$, and the limits $x_k\to x^*$ and $\rho_k\to\infty$ as $k\to\infty$, under \refa{ass21}, $B_{kj}$ is positive definite for all sufficiently large $k>0$ and $j>0$.

It follows from the Implicit Function Theorem (for example, see p.585 of \cite{NocWri99}) that there exists a $\hat\rho>0$ such that the equation
$\frac{1}{\rho_k}\na_x F(x,s_0;\mu_0,\rho_k)=0$ has a unique solution $x^*(\rho_k)$ for all $\rho_k\ge\hat\rho$, and there holds
\bea \|x^*(\rho_k)-x^*\|\le \frac{1}{\rho_k}M<\epsilon, \label{20140415e1}\eea where $\epsilon>0$ is small enough and
$$M=\max_{\|x-x^*\|<\epsilon}\left\|[\na_x^2 F(x,s_0;\mu_0,\rho_k)]^{-1}\frac{\partial}{\partial\rho}\na_x F(x,s_0;\mu_0,\rho_k)\right\|$$ is a constant independent of $\rho_k$.
\bth Under Assumptions \ref{ass2} and \ref{ass21}, if all subproblems \reff{subpro} are solved by the Newton's method, that is, at any iterate $x_{kj}$, the new iterate is generated by \bea x_{k,j+1}=x_{kj}-\frac{1}{\rho_k}B_{kj}^{-1}\na_x F(x_{kj},s_0;\mu_0,\rho_k), \label{nmd1}\eea
where $B_{kj}=\na_x^2F(x_{kj},s_0;\mu_0,\rho_k)$. Then \bea
\|x_{k+1}-x^*\|=O(\frac{1}{\rho_k})+O(\|x_{k}-x^*\|^2). \eea
Therefore, if $\frac{1}{\rho_k}=O(\|x_k-x^*\|)^2$, then the convergence is quadratic; otherwise,
if instead $\frac{1}{\rho_k}=o(\|x_k-x^*\|)$, the convergence is superlinear.
\eth\prf
Let $Q_{\rho_k}(x_{kj})=x_{kj}-\frac{1}{\rho_k}B_{kj}^{-1}\na_x F(x_{kj},s_0;\mu_0,\rho_k)$, which means that the right-hand-side is the value of function $Q_{\rho}(x)$ at $x=x_{kj}$ and $\rho=\rho_k$. Then $Q_{\rho_k}(x^*(\rho_k))=x^*(\rho_k)$, and $Q^*
=x^*$ is a limit of $Q_{\rho_k}(x_k)$ as $k\to\infty$. Note that $Q_{\rho}(x)$ is Lipschitz continuous on $\{x|\|x-x^*\|<\epsilon\}$ for any given $\epsilon>0$ and any $\rho>0$. Thus, \bea
\|x_{k,j+1}-x^*(\rho_k)\|=\|Q_{\rho_k}(x_{kj})-Q_{\rho_k}(x^*(\rho_k))\|=O(\|x_{kj}-x^*(\rho_k)\|^2), \nn
\eea
where the last equality is obtained due to $Q_{\rho_k}'(x^*(\rho_k))=0$. Therefore, by \reff{20140415e1},
\bea \dd\dd\|x_{k+1}-x^*\| \nn\\
\dd\dd=\|x_{k,j+1}-x^*(\rho_k)+x^*(\rho_k)-x^*\| \nn\\
\dd\dd\le O(\|x_{kj}-x^*(\rho_k)\|^2)+\|x^*(\rho_k)-x^*\| \nn\\
\dd\dd=O(\|x_{kj}-x^*\|^2)+O(\frac{1}{\rho_k}) \nn\\
\dd\dd=O(\|x_{k}-x^*\|^2)+O(\frac{1}{\rho_k}), \nn \eea
which completes the proof. \eop

\subsection{Linear convergence to the KKT point.}
Now we analyze the local convergence to the KKT point.
In addition to \refa{ass2}, we also need the following general conditions.
\bas\label{ass22} The following conditions hold: \\
(1) $s_k\to s^*$ and $(x^*,s^*)$ is a KKT pair; \\
(2) $s^*+c(x^*)>0$; \\
(3) matrix $B^*=\na^2f(x^*)-\sum_{i=1}^ms_i^*\na^2 c_i(x^*)+\rho_0\sum_{i\in\{i|s_i^*>0\}}\na c_i(x^*)\na c_i(x^*)^T$ is positive definite.   \eas

\refa{ass22} (1) implies $\mu_k\to 0$ as $k\to\infty$. Furthermore, by Assumptions \ref{ass2} (2) and \ref{ass22} (1), $\na^2_x F(x,s_k;\mu_k,\rho_0)$ is Lipschitz continuous. \refa{ass22} (2) is often used but is possible to be replaced by some milder condition, which is not our focus in this paper.

Our analysis in this subsection is similar to that in Bertsekas \cite{bertse82}, where he proved the linear convergence of the Hestenes-Powell augmented Lagrangian method of multipliers for nonlinear optimization with general equality constraints. We have the similar results for optimization with general inequality constraints.
\bth\label{th00} Under Assumptions \ref{ass2} and \ref{ass22}, there exist positive scalars $\bar\mu$, $\epsilon$ and $\delta$ such that, for all $(s;\mu)\in D=\{(s;\mu)|\|s-s^*\|\le\delta, \mu\in[0,\bar\mu)\}$, problem \bea \min_{x\in\{x|\|x-x^*\|<\epsilon\}}F(x,s;\mu,\rho_0) \eea
has a unique solution $x(s;\mu)$. The function $x(.;.): \Re^m\times\Re_+\to\Re^n$ is continuously differentiable in the interior of $D$, and, for all $(s;\mu)\in D$, we have \bea \|x(s;\mu)-x^*\|=O(\|s-s^*\|)+O(\|\mu\|). \label{210319a}\eea
Furthermore, if $\tilde s(s;\mu)=\rho y(x(s;\mu),s;\mu,\rho_0)$, where $y$ is a vector function with its components defined by \reff{zydf}, then there exist positive scalars $\bar\mu$, $\epsilon$ and $\delta$ such that, for all $(s;\mu)\in D=\{(s;\mu)|\|s-s^*\|\le\delta, \mu\in[0,\bar\mu)\}$, \bea
\|\tilde s(s;\mu)-s^*\|=O(\|s-s^*\|)+O(\|\mu\|). \label{210319b}\eea\eth\prf
For $\mu>0$, consider the system of equations with respect to $(x,\tilde s,s)$ and $\mu$, \bea
\dd\na f(x)-\rho_0\na c(x) y(x,s;\mu,\rho_0)=0,\dd \nn\\
\dd\tilde s-\rho_0 y(x,s;\mu,\rho_0)=0. \dd\nn\eea
By introducing the variables $t\in\Re^m$ defined by $t=s-s^*$, the preceding system can be written as \bea
\dd\na f(x)-\rho_0\na c(x) \bar y(x,t;\mu,\rho_0)=0,\dd \label{210318a}\\
\dd\tilde s-\rho_0 \bar y(x,t;\mu,\rho_0)=0, \dd\label{210318b}\eea
where $\bar y(x,t;\mu,\rho_0)=y(x,s;\mu,\rho_0)$. For $t=0$ and $\mu=0$, the system \reff{210318a}--\reff{210318b} has the solution $x=x^*$ and $\tilde s=s^*$. The Jacobian with respect to $(x,\tilde s)$ at $(x^*,s^*)$ is the order $(n+m)$ square matrix \bea
\left[\ba{cc} B^* & 0 \\
              \rho_0\hat I\na c(x^*)^T & I
\ea\right],\eea
where $I$ is the $m\times m$ identity matrix, $\hat I$ is an $m\times m$ diagonal matrix with $1$ as the $i$th diagonal element if $s_i^*>0$ and $0$ otherwise. Apparently, due to \refa{ass22} (3), it is invertible.

We now apply the implicit function theorem. It follows that there exist positive scalars $\bar\mu$, $\epsilon$ and $\delta$ and unique continuously differentiable functions $\hat x(t;\mu)$ and $\hat s(t;\mu)$ defined on $D=\{(t;\mu)|\|t\|\le\delta, \mu\in[0,\bar\mu)\}$ such that $
\|(\hat x(t;\mu)-x^*, \hat s(t;\mu)-s^*)\|<\epsilon$ for all $(t;\mu)\in D$ and satisfying \bea
\dd\na f(\hat x(t;\mu))-\na c(\hat x(t;\mu))\hat s(t;\mu)=0,\dd \label{210318c}\\
\dd\hat s(t;\mu)-\rho_0 \bar y(\hat x(t;\mu),t;\mu,\rho_0)=0. \dd\label{210318d}\eea

We differentiate \reff{210318c} and \reff{210318d} with respect to $(t;\mu)$. One obtains \bea
\left[\ba{cc} \na_t\hat x(t;\mu)^T & \hat x'(t;\mu)^T \\
              \na_t\hat s(t;\mu)^T & \hat s'(t;\mu)^T\ea\right]=B(t;\mu)\left[\ba{cc}
              0 & 0 \\
              (\hat Z(t;\mu)+\hat Y(t;\mu))^{-1}\hat Y(t;\mu) & (\hat Z(t;\mu)+\hat Y(t;\mu))^{-1} e\ea\right], \nn \eea
where $\hat x'(t;\mu)=\frac{\partial\hat x}{\partial\mu}(t;\mu)$, $\hat s'(t;\mu)=\frac{\partial\hat s}{\partial\mu}(t;\mu)$,
$\hat Z(t;\mu)=\diag(z(\hat x(t;\mu),\hat s(t;\mu);\mu,\rho_0))$, $\hat Y(t;\mu)=\diag(y(\hat x(t;\mu),\hat s(t;\mu);\mu,\rho_0))$,
\bea B(t;\mu)=\left[\ba{cc}
\na_{xx}^2L(\hat x(t;\mu),\hat s(t;\mu)) & -\na c(\hat x(t;\mu)) \\
\rho_0(\hat Z(t;\mu)+\hat Y(t;\mu))^{-1}\hat Y(t;\mu)\na c(\hat x(t;\mu))^T & I
\ea\right]^{-1}. \eea

For all $(t;\mu)$ such that $\|t\|\le\delta$ and $\mu\in [0,\bar\mu)$, we have \bea
\dd\dd\left[\ba{c} \hat x(t;\mu)-x^* \\
             \hat s(t;\mu)-s^* \ea\right]=\left[\ba{c} \hat x(t;\mu)-\hat x(0;0) \\
             \hat s(t;\mu)-\hat s(0;0) \ea\right] \nn\\
\dd\dd=\int_{0}^{1}B(\xi t;\xi\mu)\left[\ba{cc}
              0 & 0 \\
              (\hat Z(\xi t;\xi\mu)+\hat Y(\xi t;\xi \mu))^{-1}\hat Y(\xi t;\xi\mu) & (\hat Z(\xi t;\xi\mu)+\hat Y(\xi t;\xi\mu))^{-1} e\ea\right]\left[\ba{c} t\\
              \mu\ea\right]d\xi.\nn\eea
Note that $B(t;\mu)$, $(\hat Z(\xi t;\xi\mu)+\hat Y(\xi t;\xi \mu))^{-1}\hat Y(\xi t;\xi\mu)$ and $(\hat Z(\xi t;\xi\mu)+\hat Y(\xi t;\xi\mu))^{-1}$ are uniformly bounded on $D$ for all $\xi\in[0,1]$, one has \bea
\|(\hat x(t;\mu)-x^*,\hat s(t;\mu)-s^*)\|=O(\|t\|)+O(\mu), \nn\eea
which imply \reff{210319a} and \reff{210319b} immediately. \eop

We conclude this subsection by the final results, which are straightforward from \reft{th00}.
\bth Suppose that Assumptions \ref{ass2} and \ref{ass22} hold, and $\mu_k=O(\|(x_k,s_k)-(x^*,s^*)\|)$. Then \bea
{\left\|\left[\ba{c} x_{k+1}-x^*\\
s_{k+1}-s^*\ea\right]\right\|}=O{\left(\left\|\left[\ba{c} x_k-x^*\\
s_k-s^*\ea\right]\right\|\right)}. \nn\eea
Thus, the rate of convergence of the primal-dual sequence is linear. \eth

\sect{Numerical experiments}

In this paper, we mainly focus on the nonlinear programs with general nonlinear inequality constraints.
Our method can easily be extended to cope with nonlinear programs with general nonlinear inequality and equality constraints \bea
\min\dd\dd f(x) \label{prob1a-1}\\
\st\dd\dd c(x)\ge 0, \label{prob1a-2}\\
\dd\dd h(x)=0, \label{prob1a-3} \eea
with substituting for the unconstrained optimization subproblem \reff{subpro} by the optimization problem with equality constraints \bea
{\min}_x\dd\dd \frac{1}{\rho_k}F(x,s_k;\mu_k,\rho_{k}) \label{subpro-a}\\
\st\dd\dd h(x)=0. \label{subpro-c}\eea
Problem \reff{subpro-a}--\reff{subpro-c} can be approximately solved by the well developed sequential quadratic programming (SQP) methods.

The other alternative on extension of our method for program \reff{prob1a-1}--\reff{prob1a-3} is to combine with the classic augmented Lagrangian method for optimization with equality constraints, which results in the following unconstrained optimization problem \bea
{\min}_x\dd\dd \frac{1}{\rho_k}[F(x,s_k;\mu_k,\rho_{k})+\lambda_k^Th(x)+\hf\rho_k\|h(x)\|^2], \label{subpro-d}\eea
where $\lambda_k$ is an estimate of the vector of Lagrange multipliers associated with the equality constraints, and can be updated together with $s_k$ depending on the approximate solution $x_{k+1}$ of the unconstrained optimization \reff{subpro-d} by \bea
\lambda_{k+1}=\lambda_k+\rho_kh(x_{k+1}). \eea

Our targets in this section are to show that our method is usable and to demonstrate that our theoretical results are achievable.
In this sense, we will not attempt to compare our method with any recognized software, but use our method to solve some small benchmark test examples in the literature for nonlinear programs, for example, \cite{ByrCuN10,ByrMaN01,DLS17,HocSch81,LiuDaH20,LiuSun01,WacBie00}.

We have solved five examples in our numerical experiments, where the first three examples are infeasible and are solved by the proposed methods in Byrd, Curtis and Nocedal \cite{ByrCuN10} (several problems are also solved by \cite{DLS17}) for observing the rapid detection of infeasibility. The fourth and the fifth are feasible examples, where the former is challenging since the classic SQP and IPM can fail in solving the problem when starting from some infeasible points, and the latter is one for which the minimizer is a singular stationary point and the linear independence constraint qualification (LICQ) does not hold.

In our implementation, we use the standard starting points for all test problems. The initial estimates for all Lagrange multipliers (that is, the components of $s_0$) are set to be one, $H_0$ is simply set to be the identity matrix, $\mu_0=0.1$, $\rho_0=1$, $\epsilon=10^{-8}$. The subproblem \reff{subpro} is solved by the basic line search quasi-Newton method, where the step-size is computed by the Armijo line search procedure. The local convergence analysis has shown that the update of parameters can affect the rate of convergence of our method. We update the parameters $\mu_k$ and $\rho_k$ by \bea
\mu_{k+1}=\min\{0.1\mu_k, \max(\mu_k^2, \|\na F(x_{k+1},s_k;\mu_k,\rho_k)\|_{\infty}^2)\} \nn\eea
and $\rho_{k+1}=\max(\rho_k,\|s_{k+1}\|_{\infty})$ when $\tilde E_{k+1}\le 0.95\mu_k$; otherwise set $\mu_{k+1}=\mu_k$ and \bea \rho_{k+1}=\max\{2\rho_k, \min(\rho_k^2, \rho_k^2/\|\na F(x_{k+1},s_k;\mu_k,\rho_k)\|_{\infty}^2)\}. \nn \eea

The first test problem is referred as {\sl unique} in \cite{ByrCuN10}: \bea
\min\dd\dd x_1+x_2 \nn\\
({\rm TP1}) \quad\quad\st\dd\dd x_2-x_1^2-1\ge 0, \nn\\
\dd\dd 0.3(1-e^{x_2})\ge 0. \nn\eea
The standard initial point is $x_0=(3,2)$, an infeasible stationary point approximate to $x^*=(0,1)$, which is also an approximate strict minimizer of the $\ell_1$ norm of constraint violations of problem (TP1), was found in \cite{ByrCuN10}. Our algorithm terminates at $x_{11}=(0.0,0.7728)$, an approximate strict minimizer of the $\ell_2$ norm of constraint violations at which \refa{ass21} holds. The output of our algorithm is given in Table \ref{tab1}, where ``iter-sb" represents the number of iterations for  solving the $k$-th subproblem. \refal{alg1} takes the full step at all iterates including those four iterates for solving the first subproblem. It is easy to observe from Table \ref{tab1} that the convergence of $E_{k4}$ is superlinear at iterates after $k=8$ before termination. It should be noted that our algorithm only needs $14$ iterations (in every iteration we need to compute a quasi-Newton direction which is equivalent to minimizing a convex quadratic function) to find an infeasible stationary point, which is
totally different from the method presented in \cite{ByrCuN10}, where it needed to solve $24$ quadratic programming subproblems for obtaining their solution.  
\begin{table}
{\small
\begin{center}
\caption{Output for test problem (TP1): totally $14$ iterations needed.}\label{tab1} \vskip 0.2cm
\begin{tabular}{|c|c|c|c|c|c|c|c|c|}
\hline
$k$ & $f_k$ & $E_{k1}$ & $E_{k2}$ & $E_{k3}$ & $E_{k4}$ & $\mu_k$ & $\rho_k$ & iter-sb  \\
\hline
\hline
0 & 5 & 7 & 8 & 8 & 48 & 0.1 & 1 & - \\
1 & 0.1840 & 0.2381 & 0.3076 & 0.6152 & 0.5313 & 0.1 & 2 & 4 \\
2 & 0.3590 & 0.1407 & 0.1110 & 0.4440 & 0.2962 & 0.1 & 4 & 1 \\
3 & 0.4775 & 0.0385 & 0.0212 & 0.3390 & 0.1642 & 0.1 & 16 & 1 \\
4 & 0.5863 & 0.0030 & 0.0016 & 0.3533 & 0.0993 & 0.1 & 217.4755 & 1 \\
5 & 0.6776 & 0.0017 & 7.7973e-04 & 0.3601 & 0.0497 & 0.1 & 461.7971 & 1 \\
6 & 0.7602 & 7.6018e-06 & 2.7071e-06 & 0.3485 & 0.0049 & 0.1 & 1.2874e+05 & 1 \\
7 & 0.7718 & 6.9130e-07 & 2.4178e-07 & 0.3494 & 8.3494e-04 & 0.1 & 1.4451e+06 & 1 \\
8 & 0.7727 & 6.7906e-09 & 2.3745e-09 & 0.3497 & 8.2677e-05 & 0.1 & 1.4727e+08 & 1 \\
9 & 0.7728 & 3.5449e-11 & 1.2397e-11 & 0.3497 & 5.9565e-06 & 0.1 & 2.8210e+10 & 1 \\
10& 0.7728 & 9.8514e-15 & 3.4453e-15 & 0.3497 & 9.9219e-08 & 0.1 & 1.0151e+14 & 1 \\
11& 0.7728 & 1.1549e-20 & 4.0392e-21 & 0.3497 & 1.0746e-10 & 0.1 & 8.6584e+19 & 1 \\

\hline
\end{tabular}
\end{center}}
\end{table}

The second problem is the {\sl{isolated}} problem of \cite{ByrCuN10}:
\bea
\min\dd\dd x\sb{1}+x_2 \nonumber\\
({\rm TP2}) \quad\quad
\st\dd\dd -x_{1}^{2}+x\sb{2}-1\ge 0, \nonumber\\
   \dd\dd -x_{1}^{2}-x\sb{2}-1\ge 0, \nonumber\\
   \dd\dd x_{1}-x\sb{2}^2-1\ge 0, \nonumber\\
   \dd\dd -x_{1}-x\sb{2}^2-1\ge 0. \nonumber
\eea

Starting from the same initial point as problem (TP1), the algorithm in \cite{ByrCuN10} found an approximate infeasible stationary point close to $x^*=(0,0)$, a strict minimizer of the infeasibility measure in $\ell_1$ and $\ell_2$ norms. Our algorithm terminates at an approximate point $x_{10}=(4.0566, -2.4836)\times 10^{-12}$ to it. The output of our algorithm for problem (TP2) is reported in Table \ref{tab2}. The results show that the rapid convergence to the infeasible stationary point (i.e., the measurement $E_{k4}$) emerges after $k=8$. Moreover, all step-sizes are one except the $5$-th iterate for the first subproblem at which the step-size is $0.5$. We note that \cite{ByrCuN10} needs to solve $20$ quadratic programming subproblems to get the solution but our algorithm needs only $15$ iterations totally.
\begin{table}
{\small
\begin{center}
\caption{Output for test problem (TP2): totally $15$ iterations needed.}\label{tab2} \vskip 0.2cm
\begin{tabular}{|c|c|c|c|c|c|c|c|c|}
\hline
$k$ & $f_k$ & $E_{k1}$ & $E_{k2}$ & $E_{k3}$ & $E_{k4}$ & $\mu_k$ & $\rho_k$ & iter-sb  \\
\hline
\hline
0 & 5 & 13 & 12 & 12 & 126 & 0.1 & 1 & - \\
1 & -0.1912 & 0.3140 & 0.5534 & 1.1068 & 0.5927 & 0.1 & 2 & 6  \\
2 & -0.1371 & 0.1819 & 0.2684 & 1.0736 & 0.4152 & 0.1 & 4 & 1  \\
3 & -0.0797 & 0.0526 & 0.0651 & 1.0417 & 0.2412 & 0.1 & 16 & 1 \\
4 & -0.0269 & 0.0037 & 0.0040 & 1.0139 & 0.0826 & 0.1 & 256 & 1 \\
5 & -0.0023 & 1.5836e-05 & 1.5922e-05 & 1.0013 & 0.0079 & 0.1 & 6.2887e+04 & 1 \\
6 & -1.0415e-05 & 9.6070e-08 & 9.6059e-08 & 1.0001 & 3.2582e-04 & 0.1 & 1.0411e+07 & 1 \\
7 & 2.9843e-06 & 1.6318e-08 & 1.6317e-08 & 1.0000 & 1.2764e-04 & 0.1 & 6.1288e+07 & 1 \\
8 & 3.4031e-07 & 1.6238e-12 & 1.6238e-12 & 1.0000 & 1.2580e-06 & 0.1 & 6.1583e+11 & 1 \\
9 & 9.2397e-09 & 8.7481e-16 & 8.7481e-16 & 1.0000 & 2.9576e-08 & 0.1 & 1.1431e+15 & 1 \\
10& 1.5730e-12 & 5.9244e-22 & 5.9244e-22 & 1.0000 & 2.4339e-11 & 0.1 & 1.6879e+21 & 1 \\
\hline
\end{tabular}
\end{center}}
\end{table}

The third test problem is the {\sl nactive} problem in \cite{ByrCuN10}:
\begin{eqnarray}
\hbox{min} && x\sb{1} \nonumber\\
({\rm TP3}) \quad\quad
\hbox{s.t.} && \frac{1}{2}(-x\sb{1}-x\sb{2}^2-1)\ge 0, \nonumber\\
   && x\sb{1}-x_2^2\ge 0, \nonumber\\
   && -x\sb{1}+x_2^2\ge 0. \nonumber
\end{eqnarray}
This problem is still infeasible. The given initial point is $x_0=(-20,10)$. The point $x^*=(0,0)$ was an infeasible stationary point with $\|\max(0,-c(x^*))\|=0.5$. \refal{alg1} terminates at an approximate infeasible stationary point $x_8=(-0.2000,0.0000)$ with $\|\max(0,-c(\tilde x^*))\|=0.4472$, which is also a strict minimizer of the $\ell_2$ norm of the constraint violations. \refal{alg1} takes full steps at all iterates including those iterates for solving the first subproblem. It is the same as that for problems (TP1) and (TP2), the rapid convergence to the infeasible stationary point can be observed after $k=6$ from the reduction of the measurement $E_{k4}$ shown in Table \ref{tab3}. Note that \cite{ByrCuN10} needs to solve $17$ quadratic programming subproblems to get their solution.
\begin{table}
{\small
\begin{center}
\caption{Output for test problem (TP3): totally $15$ iterations needed.}\label{tab3} \vskip 0.2cm
\begin{tabular}{|c|c|c|c|c|c|c|c|c|}
\hline
$k$ & $f_k$ & $E_{k1}$ & $E_{k2}$ & $E_{k3}$ & $E_{k4}$ & $\mu_k$ & $\rho_k$ & iter-sb  \\
\hline
\hline
0 & -20 & 10 & 120 & 120 & 2805 & 0.1 & 1 & - \\
1 & -0.8948 & 0.7500 & 0.4474 & 0.8948 & 0.8685 & 0.1 & 2 & 8 \\
2 & -0.5459 & 0.3750 & 0.1365 & 0.5460 & 0.4324 & 0.1 & 4 & 1 \\
3 & -0.3673 & 0.0938 & 0.0230 & 0.3674 & 0.2092 & 0.1 & 16 & 1 \\
4 & -0.2376 & 0.0059 & 0.0015 & 0.3812 & 0.0470 & 0.1 & 256 & 1 \\
5 & -0.2020 & 2.4315e-05 & 6.4681e-06 & 0.3990 & 0.0040 & 0.1 & 6.1691e+04 & 1 \\
6 & -0.2000 & 3.8842e-06 & 1.0359e-06 & 0.4000 & 0.0016 & 0.1 & 3.8618e+05 & 1 \\
7 & -0.2000 & 1.0058e-11 & 2.6822e-12 & 0.4000 & 3.6002e-06 & 0.1 & 1.4913e+11 & 1 \\
8 & -0.2000 & 6.6858e-18 & 1.7829e-18 & 0.4000 & 2.1154e-09 & 0.1 & 2.2435e+17 & 1 \\
\hline
\end{tabular}
\end{center}}
\end{table}

In order to examine the convergence of our algorithm to the KKT point, we solve a standard test problem taken from the literature \cite{ByrMaN01,WacBie00}.
It is a well-posed problem with a single variable:
\begin{eqnarray}
\hbox{min} && x \nonumber\\
({\rm TP4}) \quad\quad
\hbox{s.t.} && x\sp{2}-1\ge 0,\nn\\
&& x-2\ge 0. \nonumber
\end{eqnarray}
This problem has a unique global minimizer $x^*=2$, at which both LICQ and MFCQ hold, and the second-order sufficient optimality conditions are satisfied. However, it is a challenging problem since a lot of infeasible starting points, such as $x_0=-4$, may bring about difficulties to the convergence to a feasible point when the problem is solved by some existing SQP and interior-point methods, for example, the linear constraints of quadratic programming subproblem of a classic SQP method at $x_0$ are inconsistent. The trouble for interior-point methods may consult \cite{LiuSun01,WacBie00}.

Starting from $x_0=-4$, our algorithm terminates at the approximate solution $x_9=2.0000$ together with $s_9=(3.3324\times 10^{-9},1.0000)$ in total $17$ iterations (including all iterations for solving subproblems). The output of \refal{alg1} is given in Table \ref{tab4}. The numbers of evaluations of functions and gradients are $18$ and all steps are full steps. The numerical results in Table \ref{tab4} show that the convergence can be much faster than the linear rate we obtained theoretically in the preceding section, but it is not so quick as that we have had for our interior-point relaxation method in \cite{LiuDaH20}, where we needed only $4$ iterations to reach a solution with higher accuracy.
\begin{table}
{\small
\begin{center}
\caption{Output for test problem (TP4): totally $17$ iterations needed.}\label{tab4} \vskip 0.2cm
\begin{tabular}{|c|c|c|c|c|c|c|c|c|}
\hline
$k$ & $f_k$ & $E_{k1}$ & $E_{k2}$ & $E_{k3}$ & $E_{k4}$ & $\mu_k$ & $\rho_k$ & iter-sb  \\
\hline
\hline
0 & -4 & 8 & 15 & 6 & 6 & 0.1 & 1 & - \\
1 & 2.2462 & 2.2462 & 2.0227 & 0 & 0 & 0.1 & 2 & 9 \\
2 & 2.1822 & 1.7978e-04 & 0.0789 & 0 & 0 & 0.01 & 2 & 1 \\
3 & 1.9439 & 0.0037 & 0.0212 & 0.0561 & 0.0561 & 0.01 & 4 & 1  \\
4 & 1.9799 & 7.9393e-04 & 0.0056 & 0.0201 & 0.0201 & 0.01 & 16 & 1  \\
5 & 2.0020 & 5.3964e-04 & 6.2535e-04 & 0 & 0 & 1.0000e-04 & 16 & 1  \\
6 & 1.9998 & 3.3671e-05 & 3.8971e-05 & 1.8455e-04 & 1.8455e-04 & 1.0000e-04 & 256 & 1 \\
7 & 2.0001 & 3.8140e-07 & 3.9063e-07 & 0 & 0 & 1.0000e-08 & 256 & 1 \\
8 & 2.0000 & 1.4898e-09 & 1.5257e-09 & 1.2935e-07 & 1.2935e-07 & 1.0000e-08 & 65536 & 1 \\
9 & 2.0000 & 8.1121e-13 & 1.5254e-13 & 0 & 0 & 2.8264e-15 & 65536 & 1 \\
\hline
\end{tabular}
\end{center}}
\end{table}

In the final experiment, we solve a standard test problem taken from \cite[Problem 13]{HocSch81}:
\begin{eqnarray}
\hbox{min} && (x\sb{1}-2)\sp{2}+x\sb{2}\sp{2} \nonumber\\
({\rm TP5}) \quad\quad
\hbox{s.t.} && (1-x\sb{1})\sp{3}-x\sb{2}\geq 0, \nonumber\\
   && x\sb{1}\geq 0, \nn\\
   && x\sb{2}\geq 0. \nonumber
\end{eqnarray}
This problem is obviously feasible, and has the optimal solution $x^*=(1,0)$ which is not a KKT point but is a singular stationary point, at which the gradients of active constraints are linearly dependent. It is a challenging problem since the convergence of many algorithms depends on the LICQ. There is not much detail on the solution of this problem in the literature.

The standard initial point $x_0=(-2,-2)$ is an infeasible point. \refal{alg1} terminates at $x_{21}\approx(1.0028,-1.0821\times 10^{-8})$, an approximate solution to the minimizer $x^*$. The parameters $\mu_k\to 0$ and $\rho_k\to\infty$ comply with the global convergence result. The output of the algorithm is given in Table \ref{tab5}. Except for the first two iterations in solving the first subproblem where the step-sizes are half, \refal{alg1} takes full steps at all other iterations.
\begin{table}
{\small
\begin{center}
\caption{Output for test problem (TP5): totally $31$ iterations needed.}\label{tab5} \vskip 0.2cm
\begin{tabular}{|c|c|c|c|c|c|c|c|c|}
\hline
$k$ & $f_k$ & $E_{k1}$ & $E_{k2}$ & $E_{k3}$ & $E_{k4}$ & $\mu_k$ & $\rho_k$ & iter-sb  \\
\hline
\hline
0 & 20 & 18 & 29 & 2 & 2 & 0.1 & 1 & - \\
1 & 0.2210 & 0.5430 & 0.7661 & 0.1048 & 0.0890 & 0.1 & 2 & 4 \\
2 & 0.2306 & 0.2847 & 0.3805 & 0.0962 & 0.0787 & 0.1 & 4 & 1 \\
3 & 0.2501 & 0.0774 & 0.0939 & 0.0790 & 0.0598 & 0.1 & 16 & 1 \\
4 & 0.2814 & 0.0054 & 0.0057 & 0.0563 & 0.0376 & 0.1 & 256 & 1 \\
5 & 0.3273 & 0.0214 & 0.0024 & 0.0443 & 0.0194 & 0.01 & 256 & 1 \\
6 & 0.5071 & 1.7996e-04 & 1.4564e-05 & 0.0120 & 0.0030 & 0.01 & 1.1631e+04 & 1 \\
7 & 0.5328 & 1.2876e-05 & 8.0585e-07 & 0.0108 & 0.0024 & 0.01 & 1.6208e+05 & 1 \\
8 & 0.5894 & 0.0040 & 6.9342e-05 & 0.0083 & 0.0040 & 0.001 & 1.6208e+05 & 1 \\
9 & 0.8225 & 1.8159e-05 & 5.7559e-08 & 0.0023 & 0.0023 & 0.001 & 3.5844e+07 & 3 \\
10& 0.8138 & 5.5583e-04 & 5.8614e-07 & 7.5580e-04 & 5.7408e-04 & 1.0000e-04 & 3.5844e+07 & 1 \\
11& 0.8390 & 3.1717e-05 & 7.4528e-09 & 5.9451e-04 & 5.9451e-04 & 1.0000e-04 & 6.2815e+08 & 2 \\
12& 0.8577 & 1.5859e-05 & 3.6279e-09 & 2.3909e-04 & 7.5132e-05 & 1.0000e-04 & 1.2563e+09 & 1 \\
13& 0.8815 & 7.9293e-06 & 1.0526e-09 & 1.3292e-04 & 3.7783e-05 & 1.0000e-04 & 2.5126e+09 & 1 \\
14& 0.9031 & 2.2753e-05 & 6.1324e-09 & 7.6758e-05 & 3.0682e-05 & 1.0000e-05 & 2.5126e+09 & 1 \\
15& 0.9535 & 2.9265e-06 & 1.0514e-10 & 1.4306e-05 & 1.4306e-05 & 1.0000e-05 & 1.9535e+10 & 2 \\
16& 0.9569 & 1.3744e-07 & 9.4547e-11 & 6.7200e-06 & 3.0639e-06 & 1.0000e-06 & 1.9535e+10 & 1 \\
17& 0.9869 & 3.8313e-07 & 5.7507e-12 & 4.0122e-07 & 4.0122e-07 & 1.0000e-07 & 1.9535e+10 & 3 \\
18& 0.9944 & 2.1689e-12 & 1.5312e-17 & 1.9388e-07 & 1.9388e-07 & 1.0000e-07 & 3.4508e+15 & 2 \\
19& 0.9944 & 1.0464e-08 & 2.5938e-16 & 1.6065e-08 & 1.0462e-08 & 1.0000e-08 & 3.4508e+15 & 1 \\
20& 0.9944 & 1.3948e-13 & 1.2128e-21 & 1.6015e-08 & 1.0402e-08 & 1.0000e-08 & 2.5889e+20 & 1 \\
21& 0.9944 & 5.1964e-14 & 8.6454e-22 & 1.0821e-08 & 3.8074e-11 & 1.0000e-08 & 6.9490e+20 & 1 \\
\hline
\end{tabular}
\end{center}}
\end{table}

In summary, the preceding numerical results not only demonstrate our strong global convergence results on \refal{alg1} for infeasible, well-posed and degenerate nonlinear programs, but also illustrate that our algorithm is capable of rapidly detecting infeasibility of nonlinear programs and can be of at least linear convergence to a KKT point of a feasible nonlinear program.

\sect{Conclusion}

The augmented Lagrangian methods of multipliers have been playing a very important role in the development of effective numerical methods for convex and nonconvex optimization problems. We present a novel augmented Lagrangian method of multipliers for nonlinear optimization with general inequality constraints. The method shares a similar algorithmic framework with existing augmented Lagrangian method of multipliers, but the subproblem has continuous second derivatives and does not depend on any projection on the set of inequality constraints. The proposed method is proved to have strong global convergence, locally it has potential of rapidly detecting the infeasibility of the original problem and can converge to the KKT point in at least a linear rate for feasible problems. The numerical experiments on some small benchmark problems demonstrate our theoretical results.

\


\end{document}